\documentclass[12pt, twoside]{article}

\usepackage{times}
\usepackage{bm}
\usepackage{graphics}
\usepackage{theorem}
\usepackage{graphicx}
\usepackage{amsmath}
\usepackage{latexsym}
\usepackage{amssymb}
\usepackage[flushmargin]{footmisc}
\renewcommand{\thefootnote}{\fnsymbol{footnote}}

\makeatletter
\def\eqnarray{\stepcounter{equation}\let\@currentlabel=\theequation
\global\@eqnswtrue
\global\@eqcnt\z@\tabskip\@centering\let\\=\@eqncr
$$\halign to \displaywidth\bgroup\@eqnsel\hskip\@centering
  $\displaystyle\tabskip\z@{##}$&\global\@eqcnt\@ne 
  \hfil$\;{##}\;$\hfil
  &\global\@eqcnt\tw@ $\displaystyle\tabskip\z@{##}$\hfil 
   \tabskip\@centering&\llap{##}\tabskip\z@\cr}
\makeatother
\makeatletter
\@addtoreset{equation}{section}
\def\theequation{\thesection.\arabic{equation}}
\makeatother
\theorembodyfont{\itshape}

\newtheorem{thm}{Theorem}[subsection]

\newtheorem{thmthm}{Theorem}[subsection]

\newtheorem{lem}{Lemma}[subsection]

\newtheorem{app-lem}{Lemma}[subsection]
\def\theapp-lem{\roman{section}.\arabic{subsection}.\arabic{app-lem}}
\newtheorem{prop}{Proposition}[subsection]

\newtheorem{app-prop}{Proposition}[subsection]
\def\theapp-prop{\roman{section}.\arabic{subsection}.\arabic{app-prop}}
\newtheorem{app-prop-except}{Proposition}[section]
\def\theapp-prop-except{\roman{section}.\arabic{app-prop-except}}

\theorembodyfont{\rmfamily}
\newenvironment{proof}{\noindent{\it Proof.~~}}{\qed}

\newtheorem{rem}{Remark}[subsection]{\it}{\rm}

\newtheorem{app-rem}{Remark}[subsection]{\it}{\rm}
\def\theapp-rem{\roman{section}.\arabic{subsection}.\arabic{rem}}
\newtheorem{defn}{Definition}[subsection]{\bf}{\rm}

\newtheorem{app-defn}{Definition}[subsection]{\bf}{\rm}
\def\theapp-defn{\roman{section}.\arabic{subsection}.\arabic{defn}}
\newtheorem{assumpt}{Assumption}{\bf}{\it}

\def\vc#1{\mbox{\boldmath $#1$}}

\newcommand{\qed}{\hspace*{\fill}$\Box$}
\newcommand{\dm}{\displaystyle}

\newcommand{\adj}{\mathrm{adj}}
\newcommand{\Mod}{\mathrm{mod}}
\newcommand{\rd}{\mathrm{d}}

\def\simhm#1{\stackrel{#1}{\sim}}

\newcommand{\rme}{\mathrm{e}}

\newcommand{\rmt}{\mathrm{t}}
\newcommand{\rmE}{\mathrm{E}}
\newcommand{\rmT}{\mathrm{T}}
\newcommand{\rsp}{\mathrm{sp}}

\newcommand{\calL}{\mathcal{L}}
\newcommand{\calS}{\mathcal{S}}

\newcommand{\bbM}{\mathbb{M}}
\newcommand{\bbN}{\mathbb{N}}
\newcommand{\bbR}{\mathbb{R}}
\newcommand{\bbS}{\mathbb{S}}

\newcommand{\bbZ}{\mathbb{Z}}

\newcommand{\EE}{\mathsf{E}}
\newcommand{\PP}{\mathsf{P}}

\newcommand{\dd}[1]{\if#11 1\!\!1 
\else {\if#1C I\!\!\!C
\else {\if#1G I\!\!\!G 
\else {\if#1J J\!\!\!J 
\else {\if#1S S\!\!\!S
\else {\if#1Z Z\!\!\!Z
\else {\if#1Q O\!\!\!\!Q
\else I\!\!#1
\fi} 
\fi}
\fi}
\fi} 
\fi} 
\fi} 
\fi} 
\renewcommand{\labelenumi}{(\alph{enumi})}

\renewcommand{\thefootnote}{\fnsymbol{footnote}}

\makeatother

\begin{document}
\thispagestyle{empty} 

\hfill
{\small Last update date: \today}

{\Large{\bf
\begin{center}
Subexponential Asymptotics of the Stationary Distributions of
  GI/G/1-Type Markov Chains\footnote[1]{This is a revised version of
the paper published in Stochastic Models vol.~29, no.~2,
pp.~190--293, 2013. \\ In the revised version, some editorial errors are
corrected and supplements are added.}
\end{center}
}
}

\begin{center}
{
Tatsuaki Kimura%
${}^{1}
\footnote[2]{Address correspondence to Tatsuaki Kimura (kimura.tatsuaki@lab.ntt.co.jp)}$, 
Hiroyuki Masuyama${}^{2}$ 
and 
Yutaka Takahashi${}^{2}$
}

\medskip

{\footnotesize
${}^1$%
NTT Network Technology Laboratories, NTT Corporation, Tokyo 180--8585 Japan

${}^2$%
Department of Systems
Science, Graduate School of Informatics, Kyoto University, Kyoto 606--8501, Japan
}

\bigskip
\medskip

{\small
\textbf{Abstract}

\medskip

\begin{tabular}{p{0.85\textwidth}}
This paper considers the subexponential asymptotics of the
  stationary distributions of GI/G/1-type Markov chains in two cases:
  (i) the phase transition matrix in non-boundary levels is
  stochastic; and (ii) it is strictly substochastic. For the case (i), we
  present a weaker sufficient condition for the subexponential
  asymptotics than those given in the literature. As for the case
  (ii), the subexponential asymptotics has not been studied, as far as
  we know. We show that the subexponential asymptotics in the case
  (ii) is different from that in the case (i). We also study the
  locally subexponential asymptotics of the stationary distributions
  in both cases (i) and (ii).
\end{tabular}
}
\end{center}

\begin{center}
\begin{tabular}{p{0.90\textwidth}}
{\small

\medskip

}
\end{tabular}

\end{center}

\section{Introduction}\label{introduction}
\renewcommand{\thefootnote}{$\ddag$\arabic{footnote}}

This paper studies the subexponential asymptotics of the stationary
distribution of an irreducible and positive recurrent Markov chain of
GI/G/1 type \cite{Gras90}. The GI/G/1-type Markov chain includes
M/G/1- and GI/M/1-type ones as special cases and plays an important
role in studying the stationary queue-length and/or waiting-time
distributions in various Markovian queues such as continuous-time
BMAP/GI/1, BMAP/D/$c$, SMAP/MSP/$c$ queues, and discrete-time
SMAP/GI/1 queues, where BMAP, SMAP and MSP represent batch Markovian
arrival process, semi-Markovian arrival process and Markovian service
process, respectively.

Let $\{(X_n,S_n);n=0,1,\dots\}$ denote a GI/G/1-type Markov chain such
that $X_n \in \bbZ_+ := \{0, 1, 2,\dots\}$ and
\[
\begin{array}{ll}
S_n \in \bbM_0 := \{1,2,\dots,M_0\}, & \mbox{if}~X_n = 0,
\\
S_n \in \bbM := \{1,2,\dots,M\}, & \mbox{otherwise},
\end{array}
\]
where $M_0$ and $M$ are positive integers. The state space of
$\{(X_n,S_n)\}$ is given by $\bbS = (\{0\} \times \bbM_0) \cup (\bbN
\times \bbM)$, where $\bbN = \{1,2,3,\dots\}$. Further, the sub-state
spaces $\{(0,j); j \in \bbM_0\}$ and $\{(k,j); j \in \bbM\}$
($k\in\bbN$) are called level 0 and level $k$, respectively.

Let $\vc{T}$ denote the transition probability matrix of the
GI/G/1-type Markov chain $\{(X_n,S_n)\}$, which can be partitioned as
follows \cite{Gras90}:
\[
\vc{T} = 
\bordermatrix{
               & {\rm lev.\ 0}&   1  &  2       &  3       & \cdots &       
\cr
{\rm lev.\ 0}\ & \vc{B}(0) & \vc{B}(1) & \vc{B}(2) & \vc{B}(3) & \cdots
\cr
\phantom{\rm lev.\ }1 & \vc{B}(-1) & \vc{A}(0) & \vc{A}(1) & \vc{A}(2) & \cdots
\cr
\phantom{\rm lev.\ }2 & \vc{B}(-2)   & \vc{A}(-1) & \vc{A}(0) & \vc{A}(1) & \cdots
\cr
\phantom{\rm lev.\ }3 & \vc{B}(-3)   & \vc{A}(-2)   & \vc{A}(-1) & \vc{A}(0) & \cdots
\cr
\phantom{\rm lev.\ }\vdots    & \vdots   & \vdots  &  \vdots   & \vdots   & \ddots
},
\]
where $\vc{A}(k)$ ($k \in \bbZ := \{ 0,\pm 1, \pm 2,\dots \}$) is an
$M \times M$ matrix, $\vc{B}(0)$ is an $M_0 \times M_0$ matrix,
$\vc{B}(k)$ ($k\in\bbN$) is an $M_0 \times M$ matrix, and $\vc{B}(k)$
($k\in\bbZ\setminus\bbZ_+$) is an $M \times M_0$ matrix. Throughout
the paper, we assume the following, unless otherwise stated.
\begin{assumpt}\label{assu-0}
(a) $\vc{T}$ is an irreducible and positive-recurrent stochastic
  matrix; (b) $\vc{A}:= \sum_{k=-\infty}^{\infty}\vc{A}(k)$ is
  irreducible.
\end{assumpt}

Under Assumption~\ref{assu-0}, $\vc{T}$ has a unique and positive
stationary distribution (see, e.g., \cite[Chapter 3,
  Theorem~3.1]{Brem99}), which is denoted by
$\vc{x}=(x_j(k))_{(k,j)\in\bbS}$. For later use, we define
$\vc{x}(0)=(x_j(0))_{j\in\bbM_0}$ and $\vc{x}(k)=(x_j(k))_{j\in\bbM}$
for $k\in\bbN$. Further, let
$\overline{\vc{x}}(k)=\sum_{l=k+1}^{\infty}\vc{x}(l)$ for
$k\in\bbZ_+$.

Some researchers have studied the subexponential asymptotics of the
stationary distribution $\vc{x}=(\vc{x}(0),\vc{x}(1),\vc{x}(2),\dots)$
of the GI/G/1-type Markov chain (including the M/G/1-type one). The
previous studies assume that $\vc{A}$ is stochastic, though $\vc{A}$
is not stochastic in general. In fact,
\[
\lim_{k\to\infty}\vc{B}(-k) \neq \vc{O}~
\mbox{if and only if}~\vc{A}\vc{e} \neq\vc{e},
\]
where $\vc{e}$ denotes a column vector of ones with an appropriate
dimension according to the context.

We briefly review the literature related to this paper. For this
purpose, let $Y$ denote a random variable in $\bbZ_+$, and for a
while, assume that
\[
\lim_{k\to \infty}{\sum_{l=k+1}^{\infty}\vc{A}(l) \over \PP(Y > k)}
= \vc{C}_1 \ge \vc{O},
~~
\lim_{k\to \infty}{\sum_{l=k+1}^{\infty}\vc{B}(l) \over \PP(Y > k)}
= \vc{C}_2  \ge \vc{O},
\]
with $\vc{C}_1 \neq \vc{O}$ or $\vc{C}_2 \neq \vc{O}$. Asmussen and
M{\o}ller~\cite{Asmu99} consider two cases: (a) $Y$ is regularly
varying; and (b) $Y$ belongs to both the subexponential class $\calS$
(see Definition~\ref{defn-subexp}) and the maximum domain of
attraction of the Gumbel distribution (see, e.g., \cite[Section
  3.3]{Embr97}). For the two cases, they show that under some
additional conditions,
\begin{equation}
\lim_{k\to\infty}{\overline{\vc{x}}(k) \over \PP(Y_{\rm e} > k)} 
= \vc{c}_1 > \vc{0}, 
\quad Y_{\rm e} \in \calS,
\label{asymp-formula}
\end{equation}
where $Y_{\rm e}$ denotes the discrete equilibrium random variable of
$Y$, distributed with $\PP(Y_{\rm e} = k) = \PP(Y > k)/\EE[Y]$
($k\in\bbZ_+$). Note here that $Y \in \calS$ does not necessarily
imply $Y_{\rm e} \in \calS$ and vice versa (see
\cite[Remark~3.5]{Sigm99}).

Li and Zhao~\cite{Li05a} show the subexponential tail asymptotics
(\ref{asymp-formula}) under the condition that $\vc{C}_2 = \vc{O}$ and
$Y$ belongs to a subclass $\calS^{\ast}$ of $\calS$ (see
Definition~\ref{def-S^*}). Note here that $Y \in \calS^{\ast}$ implies
$Y \in \calS$ and $Y_{\rm e} \in \calS$ (see Proposition A.2 in
\cite{Masu11}). Although Li and Zhao~\cite{Li05a} derive some other
asymptotic formulae for $\{\overline{\vc{x}}(k)\}$, those formulae are
incorrect due to ``the inverse of a singular matrix" (for details, see
\cite{Masu11}).

Takine~\cite{Taki04} proves that the subexponential tail asymptotics
(\ref{asymp-formula}) holds for an M/G/1-type Markov chain, assuming
that $Y_{\rm e} \in \calS$ but not necessarily $Y \in \calS$.  Thus
Takine's result shows that $Y \in \calS$ is not a necessary condition
for the subexponential decay of $\{\overline{\vc{x}}(k)\}$. However,
Masuyama \cite{Masu11} points out that Takine's proof needs an
additional condition that the $G$-matrix is aperiodic. Further,
Masuyama \cite{Masu11} presents a weaker sufficient condition for
(\ref{asymp-formula}) than those presented in the literature
\cite{Asmu99,Li05a,Taki04}, though his result is limited to the
M/G/1-type Markov chain. Recently, Kim and Kim~\cite{Kim11} improve
Masuyama \cite{Masu11}'s sufficient condition in the case where the
$G$-matrix is periodic.

In this paper, we study the subexponential decay of the tail
probabilities $\{\overline{\vc{x}}(k)\}$ in two cases: (i) $\vc{A}$ is
stochastic (i.e., $\vc{A}\vc{e} = \vc{e}$); and (ii) $\vc{A}$ is
strictly substochastic (i.e., $\vc{A}\vc{e} \le \vc{e}, \neq \vc{e}$).
For the case (i), we generalize Masuyama~\cite{Masu11}'s and Kim and
Kim~\cite{Kim11}'s results to the GI/G/1-type Markov chain. The
obtained sufficient condition for the subexponential tail asymptotics
(\ref{asymp-formula}) is weaker than those presented in Asmussen and
M{\o}ller~\cite{Asmu99} and Li and Zhao~\cite{Li05a}. As for the case
(ii), we present a subexponential asymptotic formula such that
\[
\lim_{k\to\infty}{\overline{\vc{x}}(k) \over \PP(Y > k)} 
= \vc{c}_2 > \vc{0},
\quad Y \in \calS.
\]
It should be noted that the embedded queue length process of a
BMAP/GI/1 queue with disasters falls into the case (ii) (see, e.g.,
\cite{Shin04}).  As far as we know, the subexponential asymptotics in
the case (ii) has not been studied in the literature.  Therefore, this
paper is the first report on the subexponential asymptotics in the
case (ii).

We also study the locally subexponential asymptotics of the stationary
probabilities $\{\vc{x}(k)\}$. In the case (i) (i.e., $\vc{A}$ is
stochastic), we prove the following formula under some technical
conditions:
\[
\lim_{k\to\infty}{\vc{x}(k) \over \PP(Y_{\rm e} = k)} 
= \vc{c}_3 > \vc{0},
 \quad Y \in \calS^{\ast}.
\]
Further, in the case (ii) (i.e., $\vc{A}$ is strictly substochastic),
we assume that $Y$ is {\it locally subexponential with span one}
(i.e., $Y \in \calS_{\rm loc}(1)$; see
Definition~\ref{defn-loc-S}). We then show that
\[
\lim_{k\to\infty}{\vc{x}(k) \over \PP(Y = k)} 
= \vc{c}_4 > \vc{0},
\quad Y \in \calS_{\rm loc}(1),
\]
with some technical conditions.  For the reader's convenience,
Appendix~\ref{appendix-example} presents simple examples of the case
where the stationary distribution is locally subexponential.

The rest of this paper is divided into three
sections. Section~\ref{sec-basic-results} describes some basic results
on the GI/G/1-type Markov chain and its related Markov additive
process (MAdP). In Sections \ref{sec-subexp} and \ref{sec-loc-subexp},
we studied the subexponential tail asymptotics and locally
subexponential asymptotics, respectively, of the stationary
distribution.

\section{The GI/G/1-Type Markov Chain and Its Related Markov Additive Process}
\label{sec-basic-results}

Throughout this paper, we use the following conventions. Let $\vc{I}$
denote the identity matrix with an appropriate dimension.  For any
matrix $\vc{M}$, $[\vc{M}]_{i,j}$ represents the $(i,j)$th element of
$\vc{M}$. For any matrix sequence $\{\vc{M}(k);k\in\bbZ_+\}$, let
$\overline{\vc{M}}(k) = \sum_{l=k+1}^{\infty}\vc{M}(l)$
($k\in\bbZ_+$).  For any two matrix sequences
$\{\vc{M}(k);k\in\bbZ_+\}$ and $\{\vc{N}(k);k\in\bbZ_+\}$ such that
their products are well-defined, let $\vc{M} \ast \vc{N}(k) =
\sum_{l=0}^k \vc{M}(k-l) \vc{N}(l)$ for $k\in\bbZ_+$. Further, for any
square matrix sequence $\{\vc{M}(k);k\in\bbZ_+\}$, let $\{\vc{M}^{\ast
  n}(k);k\in\bbZ_+\}$ denote the $n$-fold convolution of
$\{\vc{M}(k)\}$ with itself, i.e., $\vc{M}^{\ast n}(k) = \sum_{l=0}^k
\vc{M}^{\ast (n-1)}(k-l) \vc{M}(l)$, where $\vc{M}^{\ast0}(0) =
\vc{I}$ and $\vc{M}^{\ast0}(k) = \vc{O}$ for $k \in \bbN$.  The
conventions for matrices are also used for vectors and scalars in an
appropriate manner. Finally, the superscript ``$\rmt$" represents the
transpose operator for vectors and matrices.

\subsection{$R$- and $G$-matrices}

In this subsection, we assume that $\vc{T}$ is irreducible and
stochastic, but do not necessarily assume the recurrence of $\vc{T}$.

We consider a censored Markov chain obtained by observing
$\{(X_n,S_n)\}$ only when it is in levels 0 through $k$ ($k \in
\bbZ_+$). Let $\vc{T}^{[k]}$ ($k \in \bbZ_+$) denote the transition
probability matrix of the censored Markov chain, which is irreducible
due to the irreducibility of the original chain.  Let
$\vc{T}^{[k]}_{\nu,\eta}$ ($\nu,\eta \in \{0,1,\dots,k\}$) denote a
submatrix of $\vc{T}^{[k]}$ such that
$[\vc{T}^{[k]}_{\nu,\eta}]_{i,j}$ represents the probability that the
censored Markov chain moves from state $(\nu,i) \in \bbS$ to $(\eta,j)
\in \bbS$ in one step.

From the block Toeplitz-like structure of $\vc{T}$, we see that
$\vc{T}^{[k]}_{k-l,k}$ and $\vc{T}^{[k]}_{k,k-l}$ are independent of
$k$ if $l \in \{0,1,\dots,k-1\}$ and $k \in \bbN$. We thus define
$\vc{\Phi}(l)$ ($l \in \bbZ$) as
\begin{eqnarray}
\vc{\Phi}(l)
&=& \vc{T}^{[k]}_{k-l,k},\qquad l \in \{0,1,\dots,k-1\},~k \in \bbN,
\nonumber
\\
\vc{\Phi}(-l)
&=& \vc{T}^{[k]}_{k,k-l},\qquad l \in \{0,1,\dots,k-1\},~k \in \bbN.
\label{eqn-Phi(-l)}
\end{eqnarray}
Note here that for any fixed $\nu \in \bbN$, $[\vc{\Phi}(0)]_{i,j}$
represents the probability of hitting state $(\nu,j)$ for the first
time before entering the levels $0,1,\dots,\nu-1$, given that it
starts with state $(\nu,i)$, i.e.,
\[
[\vc{\Phi}(0)]_{i,j}
= \PP(S_{T_{\downarrow \nu}} = j \mid X_0 = \nu, S_0 = i),
\]
where $T_{\downarrow l} = \inf\{n \in \bbN; X_n = l <
X_m~(m=1,2,\dots,n-1)\}$. Thus
$\sum_{n=0}^{\infty}(\vc{\Phi}(0))^n=(\vc{I} - \vc{\Phi}(0))^{-1}$
exists because $\vc{T}^{[k]}$ is irreducible.

\begin{prop}[Theorem~1 in \cite{Gras90}]\label{prop-Gras90}
$\{\vc{\Phi}(k);k \in \bbZ\}$ is the minimal nonnegative solution of
  the following equations.
\begin{align*}
\vc{\Phi}(k) 
&= \vc{A}(k) + \sum_{m=1}^{\infty} \vc{\Phi}(k+m)
(\vc{I} - \vc{\Phi}(0))^{-1}\vc{\Phi}(-m),
& k &\in \bbZ_+,
\nonumber
\\
\vc{\Phi}(-k) 
&=
\vc{A}(-k) + \sum_{m=1}^{\infty} \vc{\Phi}(m)
(\vc{I} - \vc{\Phi}(0))^{-1}\vc{\Phi}(-k-m),
& k &\in \bbZ_+.
\end{align*}
\end{prop}

\begin{rem}\label{rem-prop-Gras90}
The proof of Theorem~1 in \cite{Gras90} is based on induction and
probabilistic interpretation, which are valid without the recurrence
of $\vc{T}$.
\end{rem}

Let $\vc{G}$ and $\vc{G}(k)$ ($k \in \bbN$) denote
\begin{equation}
\vc{G}
= \sum_{k=1}^{\infty} \vc{G}(k),
\qquad
\vc{G}(k) 
= (\vc{I} - \vc{\Phi}(0))^{-1}\vc{\Phi}(-k),
\qquad k \in \bbN,
\label{def-G(k)}
\end{equation}
respectively. Note that for any fixed $\nu \in \bbN$,
$[\vc{G}(k)]_{i,j}$ represents the probability of hitting state
$(\nu,j)$ when the Markov chain $\{(X_n,S_n)\}$ enters the levels
$0,1,\dots,\nu+k-1$ for the first time, given that it starts with
state $(\nu + k,i)$, i.e.,
\[
[\vc{G}(k)]_{i,j}
= \PP(X_{T_{<k+\nu}} = \nu, S_{T_{<k+\nu}} = j \mid X_0 = k+\nu, S_0 = i),
\qquad k \in \bbN,
\]
where $T_{<l} = \inf\{n \in \bbN; X_n < l \le X_m~(m=1,2,\dots,n-1)\}$.

Let $\vc{L}(k)$ ($k \in \bbN$) denote
\begin{equation}
\vc{L}(k) 
= \sum_{i=1}^k
\underset{n_1 + n_2 + \dots + n_i = k }
{\sum_{ (n_1,n_2,\dots,n_i) \in \bbN^i}}
\vc{G}(n_1) \vc{G}(n_2)\cdots \vc{G}(n_i), 
\qquad k \in \bbN.
\label{eqn-D(l)}
\end{equation}
For any fixed $\nu \in \bbN$, $[\vc{L}(k)]_{i,j}$ represents the
probability of hitting state $(\nu,j)$ when the Markov chain
$\{(X_n,S_n)\}$ enters the levels $0,1,\dots,\nu$ for the first time,
given that it starts with state $(\nu + k,i)$, i.e.,
\[
[\vc{L}(k)]_{i,j}
=\PP(S_{T_{\downarrow \nu}} = j \mid X_0 = k+\nu, S_0 = i).
\]
It follows from (\ref{eqn-D(l)}) that
\begin{equation}
\widehat{\vc{L}}(z) 
:= \sum_{k=1}^{\infty} z^{-k} \vc{L}(k)
= ( \vc{I} - \widehat{\vc{G}}(z) )^{-1} \widehat{\vc{G}}(z),
\label{defn-hat{L}(z)}
\end{equation}
where $\widehat{\vc{G}}(z) = \sum_{k=1}^{\infty} z^{-k} \vc{G}(k)$.

Let $\vc{R}_0(k)$ and $\vc{R}(k)$ ($k \in \bbZ_+$) denote $M_0 \times
M$ and $M \times M$ matrices, respectively, such that
\begin{align}
\vc{R}_0(0) &= \vc{O}, & \vc{R}(0) &= \vc{O}, \nonumber
\\
\vc{R}_0(k) 
&= \vc{T}_{0,k}^{[k]}(\vc{I} - \vc{\Phi}(0))^{-1},
&
\vc{R}(k) 
&= \vc{\Phi}(k)(\vc{I} - \vc{\Phi}(0))^{-1},
\quad  k \in \bbN.
\label{def-R(k)}
\end{align}
For any fixed $\nu \in \bbN$, $[\vc{R}(k)]_{i,j}$ ($k \in \bbN$)
represents the expected number of visits to state $(\nu+k,j)$ before
entering the levels $0,1,\dots,\nu+k-1$, given that the Markov chain
$\{(X_n,S_n)\}$ starts with state $(\nu,i)$.  Further, $\vc{R}_0(k)$
($k \in \bbN$) can be interpreted in the same way though $\nu \in
\bbN$ is replaced by zero. Formally, for $k \in \bbN$,
\begin{eqnarray*}
[\vc{R}_0(k)]_{i,j}
&=& \EE\left[ \left. \sum_{n=1}^{T_{< k}} \dd{1}(X_n = k, S_n = j) 
\,\right| X_0 = 0, S_0 = i \right],
\\
\phantom{}[\vc{R}(k)]_{i,j}
&=& \EE\left[ \left. \sum_{n=1}^{T_{< k+\nu}} \dd{1}(X_n = k+\nu, S_n = j) 
\,\right| X_0 = \nu \in \bbN, S_0 = i \right],
\end{eqnarray*}
where $\dd{1}(\chi)$ denotes the indicator function of an event
$\chi$. It follows from the definitions of $\vc{R}_0(k)$, $\vc{R}(k)$,
$\vc{L}(k)$ and $\vc{\Phi}(0)$ that
\begin{eqnarray}
\vc{R}_0(k)
&=&
\left[
\vc{B}(k) + \sum_{m=1}^{\infty} \vc{B}(k+m) \vc{L}(m) \right]
(\vc{I} - \vc{\Phi}(0))^{-1},
\quad k \in \bbN,
\label{eqn-R0-B}
\\
\vc{R}(k) 
&=&
\left[
\vc{A}(k)  + \sum_{m=1}^{\infty} \vc{A}(k+m) \vc{L}(m)\right]
(\vc{I} - \vc{\Phi}(0))^{-1}, 
\quad k \in \bbN,
\label{eqn-R-A}
\end{eqnarray}
which hold without the recurrence of $\vc{T}$.
 
We now define $\widehat{\vc{R}}_0(z)$, $\widehat{\vc{R}}(z)$ and
$\widehat{\vc{B}}(z)$ as
\[
\widehat{\vc{R}}_0(z)
= \sum_{k=1}^{\infty} z^k \vc{R}_0(k),
~~
\widehat{\vc{R}}(z)
= \sum_{k=1}^{\infty} z^k \vc{R}(k),
~~
\widehat{\vc{B}}(z)
= \sum_{k=1}^{\infty} z^k \vc{B}(k),
\]
respectively. We then have the following result.
\begin{prop}[Theorem~1 and Lemma~3 in \cite{Li05b}]
\label{prop-convergence-radii}
Let $r_{R_0}$, $r_R$, $r_G$, $r_{A_+}$, $r_{A_-}$ and $r_B$ denote the
convergence radii of $\widehat{\vc{R}}_0(z)$, $\widehat{\vc{R}}(z)$,
$\widehat{\vc{G}}(1/z)=\sum_{k=1}^{\infty} z^k \vc{G}(k)$,
$\sum_{k=1}^{\infty} z^k \vc{A}(k)$, $\sum_{k=1}^{\infty} z^k
\vc{A}(-k)$ and $\widehat{\vc{B}}(z)$, respectively. Then $r_{R_0} =
r_B \ge 1$, $r_R = r_{A_+} \ge 1$ and $r_G = r_{A_-} \ge 1$.
\end{prop}

\begin{prop}[Theorem~14 in \cite{Zhao03}]\label{prop-RG}
Let $\widehat{\vc{A}}(z) = \sum_{k \in \bbZ}
z^k \vc{A}(k)$. We then have
\begin{equation}
\vc{I} - \widehat{\vc{A}}(z)
= (\vc{I} - \widehat{\vc{R}}(z)) (\vc{I} - \vc{\Phi}(0)) 
(\vc{I} - \widehat{\vc{G}}(z)),
\qquad
|z|  \in I_A,
\label{eq-RG-1}
\end{equation}
where $I_A = (1/r_{A_-},r_{A_+}) \cup \{1\}$. 
\end{prop}

\begin{rem}
Although Theorem~14 in \cite{Zhao03} assume that $\vc{A}$ is
irreducible and stochastic, these conditions are not necessarily
required by the algebraic proof of the theorem.
\end{rem}

\begin{prop}\label{prop-Phi}
Let $\vc{R} = \sum_{k=1}^{\infty}\vc{R}(k)$.  If $\vc{A}$ is
irreducible and strictly substochastic, then (i) $\rsp(\vc{G}) < 1$;
(ii) $\rsp(\vc{R}) < 1$; and (iii) $\rsp(\sum_{l=0}^{\infty}
\vc{\Phi}(-l)) < 1$, where $\rsp(\,\cdot\,)$ denotes the spectral radius
of a matrix in parentheses.
\end{prop}

\proof See Appendix~\ref{proof-prop-Phi}. \qed

\subsection{Sufficient conditions for positive recurrence}

In this subsection, we provide two sets of sufficient conditions for
Assumption~\ref{assu-0}.  For later use, let $\vc{\pi} > \vc{0}$
denote a left eigenvector of $\vc{A}$ such that $\vc{\pi}\vc{A} =
\rsp(\vc{A})\vc{\pi}$ and $\vc{\pi}\vc{e} = 1$ (see Theorem~8.4.4 in
\cite{Horn90}). Let $\sigma$ denote
\begin{equation}
\sigma =\vc{\pi}\sum_{k \in \bbZ} k\vc{A}(k)\vc{e}.
\label{defn-sigma}
\end{equation}
If $\vc {A}$ is stochastic, then $\vc{\pi}$ is the unique invariant
probability vector of $\vc{A}$ and $\sigma$ is the conditional mean
drift of the level process $\{X_n;n\in\bbZ_+\}$ with $X_n \ge 1$.

\begin{prop}[Proposition 3.1 in Chapter XI of \cite{Asmu03-book}]
\label{prop-stability-1}
Suppose $\vc{T}$ and $\vc{A}$ are irreducible and stochastic.  Then
$\vc{T}$ is positive recurrent if and only if $\sigma < 0$ and
$\sum_{k=1}^{\infty} k\vc{B}(k)\vc{e} < \infty$.
\end{prop}

\begin{prop}\label{prop-stability-2}
Suppose $\vc{T}$ is irreducible and stochastic. Then if $\vc{A}$ is
irreducible and strictly substochastic, $\vc{T}$ is positive
recurrent.
\end{prop}

\proof Proposition~\ref{prop-Phi} implies that
$\lim_{k\to\infty}\vc{R}^k = \vc{O}$ and $(\vc{I} - \vc{G})^{-1}$
exists. Further from (\ref{eqn-R0-B}), we have
\begin{eqnarray*}
\vc{R}_0
&:=& \sum_{k=1}^{\infty}\vc{R}_0(k)
\nonumber
\\
&=&
\left[
\sum_{k=1}^{\infty}\vc{B}(k) 
+ \sum_{m=1}^{\infty} \left(\sum_{k=1}^{\infty}\vc{B}(k+m) \right) \vc{L}(m) \right]
(\vc{I} - \vc{\Phi}(0))^{-1}
\nonumber
\\
&\le&
\sum_{k=1}^{\infty}\vc{B}(k)
\left[
\vc{I}  + \sum_{m=1}^{\infty} \vc{L}(m) \right]
(\vc{I} - \vc{\Phi}(0))^{-1}
\nonumber
\\
&=&
\sum_{k=1}^{\infty}\vc{B}(k)
(\vc{I} - \vc{G})^{-1}(\vc{I} - \vc{\Phi}(0))^{-1} < \infty,
\end{eqnarray*}
where the last equality follows from (\ref{defn-hat{L}(z)}). As a
result, it follows from Theorem~3.4 in \cite{Zhao98} that $\vc{T}$ is
positive recurrent.  \qed

\subsection{Matrix-product form of the stationary distribution}

This subsection discusses the stationary distribution $\{\vc{x}(k)\}$
under Assumption~\ref{assu-0}.  It is easy to see that
$(\vc{x}(0),\vc{x}(1),\dots,\vc{x}(k))$ is an invariant measure vector
of the censored transition probability matrix $\vc{T}^{[k]}$, i.e.,
\[
(\vc{x}(0),\vc{x}(1),\dots,\vc{x}(k))\vc{T}^{[k]}
= (\vc{x}(0),\vc{x}(1),\dots,\vc{x}(k)),
\]
which leads to
\begin{equation}
\vc{x}(k)
= \left[\vc{x}(0) \vc{T}^{[k]}_{0,k}
+ \sum_{l=1}^{k-1} \vc{x}(l) \vc{T}^{[k]}_{l,k}\right] 
(\vc{I} - \vc{T}^{[k]}_{k,k})^{-1},
\qquad k \in \bbN.
\label{eqn-x(k)}
\end{equation}
In terms of $\vc{R}(k)$ and $\vc{R}_0(k)$, we can rewrite
(\ref{eqn-x(k)}) as
\begin{equation}
\vc{x}(k) 
= \vc{x}(0) \vc{R}_0(k) + \sum_{l=1}^k \vc{x}(l) \vc{R}(k-l),
\qquad k \in \bbN,
\label{eq-x-R-R0}
\end{equation}
where we use $\vc{R}(0)=\vc{O}$. It then follows from
(\ref{eq-x-R-R0}) that
\begin{equation}
\vc{x}(k) = \vc{x}(0) \vc{R}_0 * \vc{F}(k),\qquad k \in \bbN,
\label{eqn-x(k)-convolution}
\end{equation}
where $\vc{F}(k)$ $(k\in\bbZ_+)$
is given by
\begin{equation}
\vc{F}(k) = \sum_{n=0}^{\infty}\vc{R}^{\ast n}(k). 
\label{def-F(k)}
\end{equation}
Thus we have
\begin{equation}
\overline{\vc{x}}(k) =
\vc{x}(0) \overline{ \vc{R}_0 * \vc{F}}(k), \qquad k\in\bbZ_+.
\label{eq-xk-R0-Gamma}
\end{equation}
Further let $\widehat{\vc{x}}(z) =
\sum_{k=1}^{\infty}z^k \vc{x}(k)$. We then have
\begin{equation}
\widehat{\vc{x}}(z) 
= \vc{x}(0) \widehat{\vc{R}}_0(z) 
( \vc{I} - \widehat{\vc{R}}(z) )^{-1}.
\label{eq-RG-2}
\end{equation}
Letting $z=1$ in (\ref{eq-RG-2}) yields
\begin{equation}
\overline{\vc{x}}(0) = 
\vc{x}(0) \vc{R}_0 (\vc{I} - \vc{R} )^{-1},
\label{eq-bar-x_0-R}
\end{equation}
where $\vc{R}_0 = \sum_{k=1}^{\infty}\vc{R}_0(k)$.

\subsection{Period of the related Markov additive process}\label{subsec-period-MAdP}

We consider a MAdP $\{(\breve{X}_n, \breve{S}_n);n\in\bbZ_+\}$ with
state space $\bbZ \times \bbM$ and kernel $\{\vc{A}(k);k\in
\bbZ\}$. The stochastic behavior of the MAdP $\{(\breve{X}_n,
\breve{S}_n)\}$ is equivalent to that of the GI/GI/1-type Markov chain
$\{(X_n,S_n)\}$ while the latter is being in non-boundary levels,
i.e., for any $i,j \in \bbM$,
\begin{eqnarray}
&& \PP(\breve{X}_{n+1}=k,~\breve{S}_{n+1} = j 
\mid \breve{X}_n=l,~\breve{S}_n = i)
\nonumber
\\
&& \quad {} 
= \PP(X_{n+1}=k,~S_{n+1} = j \mid X_n=l,~S_n = i),
\quad k,l \in \bbN.
\qquad
\label{add-eqn-58}
\end{eqnarray}
The period of the MAdP $\{(\breve{X}_n, \breve{S}_n)\}$, denoted by
$\tau$, is the largest positive integer such that
\begin{equation}
[\vc{A}(k)]_{i,j} > 0
\mbox{ only if } k \equiv
p(j)-p(i)~(\Mod~\tau),
\label{add-eqn-12}
\end{equation}
where $p$ is some function $p$ from $\bbM$ to $\{0,1,\dots,\tau-1\}$
(see Appendix~B in \cite{Kimu10} and its revised version
\cite{Kimu11}).

\begin{rem}
Lemma~B.2 in \cite{Kimu10} states that function $p$ satisfying
(\ref{add-eqn-12}) is {\it injective}, which is not true in
general. This error is corrected in the revised version \cite{Kimu11}.
\end{rem}

\begin{rem}
If the Markov chain $\{(X_n,S_n)\}$ is of M/G/1 type, the period
$\tau$ is less than or equal to $M$ (see, e.g., Proposition~2.9 in
\cite{Kimu10}).
\end{rem}

\begin{rem}
We suppose 
\begin{eqnarray*}
\vc{A}(0) &=& 
\vc{O},
\quad
\vc{A}(1) = 
\left(
\begin{array}{ccc}
0 & 0 & \dm{1 \over 6} 
\\
\rule{0mm}{7mm}0 & 0 & \dm{1 \over 6} 
\\
\rule{0mm}{7mm}\dm{1 \over 6} & \dm{1 \over 6}& 0
\end{array}
\right),
\\
\vc{A}(-2) &=& 
\left(
\begin{array}{ccc}
\dm{1 \over 3}  & \dm{1 \over 3} & 0
\\
\rule{0mm}{7mm}\dm{1 \over 3} & \dm{ 1 \over 3} & 0
\\
\rule{0mm}{7mm}0 & 0 &  \dm{1 \over 3}
\end{array} 
\right),\quad
\vc{A}(-1) = 
\left(
\begin{array}{ccc}
0 & 0 &  \dm{ 1 \over 6}
\\
\rule{0mm}{7mm}0 & 0 &  \dm{1 \over 6}
\\
\rule{0mm}{7mm}\dm{1 \over 6} & \dm{1 \over 6} & 0
\end{array}
\right).
\end{eqnarray*}
Let $p(0) = p(1) = 1$ and $p(2) = 0$. It then follows that
\[
[\vc{A}(k)]_{i,j} > 0 \mbox{ only if } k \equiv p(j) - p(i)\ (\mbox{mod } 2),
\] 
and thus the period of MAdP with kernel $\{\vc{A}(k)\}$ is equal to
two.
\end{rem}

We now introduce the following notation.
\begin{defn}\label{def-delta}
For any finite square matrix $\vc{X}$ with possibly complex elements,
let $\delta(\vc{X})$ denote an eigenvalue of $\vc{X}$, which satisfies
$\delta(\vc{X}) = \rsp(\vc{X})e^{\iota \xi}$ and
\[
\xi
= \inf\{0 \le x < 2\pi; \det(\rsp(\vc{X})e^{\iota x}\vc{I} - \vc{X}) = 0\},
\]
where $\iota$ denotes the imaginary unit, i.e., $\iota = \sqrt{-1}$.
\end{defn}

\begin{rem}
Suppose $\vc{X}$ is nonnegative. We then have $\delta(\vc{X}) =
\rsp(\vc{X})$ (see Theorem~8.3.1 in \cite{Horn90}). Further, if $\vc{X}$
is irreducible, $\delta(\vc{X})$ is the Perron-Frobenius eigenvalue of
$\vc{X}$ (see Theorem~8.4.4 in \cite{Horn90}).
\end{rem}

Let $\vc{\mu}(z)$ and $\vc{v}(z)$ denote the left- and
right-eigenvectors of $\widehat{\vc{A}}(z)$ corresponding to the
eigenvalue $\delta(\widehat{\vc{A}}(z))$, normalized such that
\[
\vc{\mu}(z) \vc{\Delta}_M(z/|z|)\vc{e} = 1, 
\qquad 
\vc{\mu}(z)\vc{v}(z) = 1,
\]
where $\vc{\Delta}_M(z)$ denotes an $M \times M$ diagonal matrix as
follows:
\[
\vc{\Delta}_M(z) = \left(
\begin{array}{cccc}
z^{-p(1)} &          &         &
\\
         & z^{-p(2)} &         &
\\
         &          &\ddots   &  
\\
         &          &         & z^{-p(M)}
\end{array}
\right).
\]
Note that $\vc{\mu}(1) = \vc{\pi}$ and $\vc{v}(1) = \vc{e}$. Further,
let $\omega$ denote an arbitrary complex number such that $|\omega| =
1$. We then have the following results.
\begin{prop}[Lemma B.3 in \cite{Kimu10}]\label{prop-delta-Gamma_A-1}
Suppose Assumption~\ref{assu-0} holds and let $\omega_x =
\exp(2\pi\iota/x)$ for $x\ge1$. Then the following are true for all $y
\in I_A$ and $\nu=0,1,\dots,\tau-1$.
\begin{enumerate}\renewcommand{\labelenumi}{(\roman{enumi})}
\item $\delta(\widehat{\vc{A}}(y\omega_{\tau}^{\nu})) =
  \delta(\widehat{\vc{A}}(y))$, both of which are simple eigenvalues;
  and
\item $\vc{\mu}(y\omega_{\tau}^{\nu}) =
  \vc{\mu}(y)\vc{\Delta}_M(\omega_{\tau}^{\nu})^{-1}$ and
  $\vc{v}(y\omega_{\tau}^{\nu}) =
  \vc{\Delta}_M(\omega_{\tau}^{\nu})\vc{v}(y)$.
\end{enumerate}
\end{prop}
\begin{prop}[Theorem B.1 in \cite{Kimu10}]\label{prop-MAdP}
Suppose Assumption~\ref{assu-0} holds and $\delta(\widehat{\vc{A}}(y))
= 1$ for some $y \in I_A$. Then $\delta(\widehat{\vc{A}}(y\omega)) =
1$ if and only if $\omega^{\tau} = 1$. Therefore
\[
\tau = \max\{n \in \bbN; \delta(\widehat{\vc{A}}(y\omega_n)) = 1\}.
\]
Further if $\delta(\widehat{\vc{A}}(y\omega)) = 1$, the eigenvalue is
simple.
\end{prop}


\subsection{Spectral analysis of $G$-matrices from stochastic $\vc{A}$}

In this subsection, we assume that Assumption~\ref{assu-0} holds and
$\vc{A}$ is stochastic. Under the assumptions, $\vc{G}$ is stochastic,
i.e., $\delta(\vc{G}) = 1$ (see Theorem 3.4 in \cite{Zhao98}).

We first provide a basic result on the structure of $\vc{G}$.
\begin{prop}\label{prop-structure-G}
Suppose Assumption~\ref{assu-0} holds and $\vc{A}$ is stochastic. Then
$\vc{G}$ has an exactly one irreducible class, denoted by
$\bbM_{\bullet} \subseteq \bbM$. Thus, $\vc{G}$ is irreducible, or
after some permutations it takes a form such that
\[
\bordermatrix{ 
       &          \bbM_{\bullet}    &   \bbM_{\rm T}
\cr  \bbM_{\bullet}  &  \vc{G}_{\bullet}  &  \vc{O} 
\cr  \bbM_{\rm T}          &  \vc{G}_{\circ}  & \vc{G}_{\rm T} 
},    
\qquad \bbM_{\rm T}:=\bbM \setminus \bbM_{\bullet},
\]
where
$\vc{G}_{\bullet}$ is irreducible, $\vc{G}_{\rm T}$ is strictly lower
triangular and $\vc{G}_{\circ}$ does not have, in general, a special
structure.
\end{prop}

\proof See Appendix~\ref{proof-prop-structure-G}.  \qed

\begin{rem}\hspace{-1.5mm}\footnote{This remark is added in the revised version.}\label{rem-G-matrix}
The proof of Proposition~\ref{prop-structure-G} relies only on the facts that (i) $\vc{A}$ is irreducible; and (ii) $\vc{G}$ is not a nilpotent matrix. 
Thus, the irreducibility of $\vc{A}$ and ${\rm sp}(\vc{G}) > 0$ imply that Proposition~\ref{prop-structure-G} holds. On the other hand, Proposition~\ref{prop-RG} yields 
\[
\det(\vc{I} - \vc{A}) 
= \det(\vc{I} - \vc{R}) \det(\vc{I} - \vc{\Phi}(0)) \det(\vc{I} - \vc{G}).
\]
Note here that if $\vc{T}$ is positive recurrent then $\det(\vc{I} - \vc{R}) \neq 0$ (see the proof of Theorem 3.4 of \cite{Zhao98}). Therefore, if the conditions of Proposition~\ref{prop-structure-G} are satisfied, then $\det(\vc{I} - \vc{A}) = 0$, $\det(\vc{I} - \vc{\Phi}(0)) \neq 0$ and $\det(\vc{I} - \vc{R}) \neq 0$ and thus
$\det(\vc{I} - \vc{G}) = 0$, which implies that ${\rm sp}(\vc{G}) = 1$.
\end{rem}

\medskip

Let $\vc{G}_{\bullet}(k)$ $(k\in\bbN)$ denote the square submatrix of
$\vc{G}(k)$ $(k\in\bbN)$ corresponding to the irreducible class
$\bbM_{\bullet} \subseteq \bbM$, i.e., $\vc{G}_{\bullet}=
\sum_{k=1}^{\infty}\vc{G}_{\bullet}(k)$. Further let
$\widehat{\vc{G}}_{\bullet}(z) =
\sum_{k=1}^{\infty}z^{-k}\vc{G}_{\bullet}(k)$. It then follows from
Proposition~\ref{prop-structure-G} that
\begin{equation}
\delta(\widehat{\vc{G}}(z)) = \delta(\widehat{\vc{G}}_{\bullet}(z)),
\label{add-eqn-04}
\end{equation}
because $\vc{G}_{\rm T}$ (if any) is a nilpotent matrix.

We now consider a MAdP $\{(\breve{X}^{(G)}_n,
\breve{S}^{(G)}_n);n\in\bbZ_+\}$ with state space $\bbZ \times
\bbM_{\bullet}$ and kernel $\{\vc{\varGamma}^{(G)}(k);k\in \bbZ\}$,
where
\begin{equation}
\vc{\varGamma}^{(G)}(k) = \left\{
\begin{array}{ll}
\vc{O}, & k \in\bbZ_+,
\\
\vc{G}_{\bullet}(k),  &  k \in  \bbZ\backslash\bbZ_+.
\end{array}
\right.
\label{def-Gamma^G(k)}
\end{equation}
Equation (\ref{def-Gamma^G(k)}) and the irreducibility of
$\sum_{k\in\bbZ}\vc{\varGamma}^{(G)}(k)=\vc{G}_{\bullet}$ imply that
the period of the MAdP $\{(\breve{X}^{(G)}_n, \breve{S}^{(G)}_n)\}$ is
well-defined (see Definition~B.1 in \cite{Kimu10}) and denoted by
$\tau_G$. Combining (\ref{add-eqn-04}), (\ref{def-Gamma^G(k)}) and
Theorem~B.1 in \cite{Kimu10}, we obtain
\begin{equation}
\tau_G = \max\{n \in \bbN; \delta(\widehat{\vc{G}}(\omega_n)) = 1\}.
\label{def-tau_G}
\end{equation}

\begin{prop}\label{prop-period-G-a}
Suppose Assumption~\ref{assu-0} holds and $\vc{A}$ is stochastic. Then
the following are true.
\begin{enumerate}
\renewcommand{\labelenumi}{(\roman{enumi})}
\item $\tau_G = \tau$;
\item $\delta(\widehat{\vc{G}}(\omega)) = 1$ if and only if
  $\omega^{\tau} = 1$;
\item if $\delta(\widehat{\vc{G}}(\omega)) = 1$,
  the eigenvalue is simple; and
\item for $y > 1/r_{A_-}$,
\[
\delta(\widehat{\vc{G}}(y\omega_{\tau}^{\nu})) 
= \delta(\widehat{\vc{G}}(y)),
\quad \nu=0,1,\dots,\tau-1,
\]
which are simple eigenvalues of
$\widehat{\vc{G}}(y\omega_{\tau}^{\nu})$ and
$\widehat{\vc{G}}(y)$, respectively.
\end{enumerate}

\end{prop}

\proof See Appendix~\ref{proof-prop-period-G-a}. \qed


\medskip

We define $\lambda_i^{(G)}(z)$'s $(i=2,3,\dots,M)$ as the eigenvalues
of $\widehat{\vc{G}}(z)$ such that $\delta(\widehat{\vc{G}}(z)) \ge
|\lambda_i^{(G)}(z)|$ (see Definition \ref{def-delta}). We then have
\begin{equation}
\det(\vc{I} - \widehat{\vc{G}}(z)) = (1 - \delta(\widehat{\vc{G}}(z)))
\prod_{i=2}^M (1 - \lambda_i^{(G)}(z)).
\label{def-eigenvalues-G}
\end{equation}

\begin{prop}\label{lem-adjG}
Suppose Assumption \ref{assu-0} holds and $\vc{A}$ is stochastic. Let
\begin{align}
\vc{\psi}(\omega_{\tau}^{\nu}) 
&= {\vc{\pi} (\vc{I} - \vc{R})( \vc{I} - \vc{\Phi}(0))
\over
\vc{\pi}(\vc{I} - \vc{R}) 
( \vc{I} - \vc{\Phi}(0))\vc{e} }\vc{\Delta}_M(\omega_{\tau}^{\nu})^{-1} ,
& \nu &= 0,1,\dots,\tau-1,
\label{eq-left-G}
\\
\vc{y}(\omega_{\tau}^{\nu}) 
&= \vc{\Delta}_M(\omega_{\tau}^{\nu})\vc{e},
& \nu &= 0,1,\dots,\tau-1.
\label{eq-right-G}
\end{align}
Then the following hold for $\nu=0,1,\dots,\tau-1$: (i)
$\vc{\psi}(\omega_{\tau}^{\nu})$ and $\vc{y}(\omega_{\tau}^{\nu})$ are
the left- and right-eigenvectors of
$\widehat{\vc{G}}(\omega_{\tau}^{\nu})$ corresponding to the eigenvalue
$\delta(\widehat{\vc{G}}(\omega_{\tau}^{\nu}))=1$; and (ii)
\[
\adj (\vc{I} - \widehat{\vc{G}}(\omega_{\tau}^{\nu} ))
=
\prod_{i=2}^M ( 1 - \lambda_i^{(G)}(\omega_{\tau}^{\nu}) ) 
\vc{y}( \omega_{\tau}^{\nu}) \vc{\psi} (\omega_{\tau}^{\nu}), 
\]
where $\adj(\vc{Y})$ denotes the
adjugate matrix of a square matrix $\vc{Y}$.
\end{prop}

\proof See Appendix~\ref{proof-lem-adjG}. \qed

\section{Subexponential Tail Asymptotics}\label{sec-subexp}

This section studies the subexponential decay of the tail
probabilities $\{\overline{\vc{x}}(k)\}$, under the following
assumption.
\begin{assumpt}\label{assu-sec3}
Either of (I) and (II) is satisfied:
\begin{enumerate}
\renewcommand{\labelenumi}{(\Roman{enumi})}
\item Assumption~\ref{assu-0} holds, $\vc{A}$ is stochastic,
  and $\sum_{k\in\bbZ}|k|\vc{A}(k) < \infty$; or
\item Assumption~\ref{assu-0} holds and $\vc{A}$ is strictly
  substochastic.
\end{enumerate}

\end{assumpt}
Assumption~\ref{assu-sec3}~(I) and (II) are considered in
subsections~\ref{subsec-stochastic} and \ref{subsec-substochastic},
respectively.

\subsection{Case of stochastic \vc{A}}\label{subsec-stochastic}

\begin{lem}\label{lem-sigma}
Under Assumption~\ref{assu-sec3}~(I), 
\begin{equation}
\sigma = -
\vc{\pi}
(\vc{I} - \vc{R} ) (\vc{I} - \vc{\Phi}(0)) 
\sum_{k=1}^{\infty}k\vc{G}(k) \vc{e} \in (-\infty,0),
\label{eqn-sigma}
\end{equation}
where $\sigma$ is defined in (\ref{defn-sigma}).
\end{lem}

\proof We have $-\infty < \sigma < 0$ due to (\ref{defn-sigma}),
Proposition~\ref{prop-stability-1} and the third condition of
Assumption~\ref{assu-sec3}~(I). Further since $\sigma = \vc{\pi}(\rd /
\rd z)\widehat{\vc{A}}(z)|_{z=1} \vc{e}$ and $(\rd / \rd
z)\widehat{\vc{G}}(z)|_{z=1}=-\sum_{k=1}^{\infty}k\vc{G}(k)$, we
obtain (\ref{eqn-sigma}) by differentiating (\ref{eq-RG-1}) with
respect to $z$, pre-multiplying by $\vc{\pi}$, post-multiplying by
$\vc{e}$ and letting $z=1$.  \qed

\medskip

Using Lemma~\ref{lem-sigma}, Propositions~\ref{prop-period-G-a} and
\ref{lem-adjG}, we obtain the following result.
\begin{lem}\label{lem-lim-D}
If Assumption~\ref{assu-sec3}~(I) holds, then for
$l=0,1,\dots,\tau-1$,
\begin{equation}
\lim_{n\to \infty} \vc{L}(n\tau + l) 
=
\sum_{\nu=0}^{\tau-1}{1 \over ( \omega_{\tau}^{\nu} )^l}
\vc{\Delta}_M(\omega_{\tau}^{-\nu} )
\vc{e}\vc{\psi}
\vc{\Delta}_M(\omega_{\tau}^{-\nu})^{-1},
\label{eq-lim-D}
\end{equation}
where 
\begin{equation}
\vc{\psi} = \vc{\pi}(\vc{I} - \vc{R} ) (\vc{I} -
\vc{\Phi}(0))/(-\sigma).
\label{defn-psi}
\end{equation}

\end{lem}
\proof See Appendix~\ref{proof-lem-lim-D}.\qed

\medskip

For $l=0,1,\dots,\tau-1$, let $\bbM^{(l)} = \{j\in\bbM;p(j) = l\}$ and
$|\bbM^{(l)}|$ denote the cardinality of $\bbM^{(l)}$. Further, let
$\vc{\psi}^{(l)}$ denote a subvector of $\vc{\psi}$ corresponding to
$\bbM^{(l)}$, and $\vc{e}^{(l)}$ denote an $|\bbM^{(l)}| \times 1$
vector of ones. Note here that
$\sum_{\nu=0}^{\tau-1}(\omega_{\tau}^m)^{\nu} = 0$ for all $m =
1,2,\dots,\tau-1$ because $\omega_{\tau},
\omega_{\tau}^2,\dots,\omega_{\tau}^{\tau-1}$ are the solutions of the
equation $\sum_{\nu=0}^{\tau-1}z^{\nu} = 0$. It then follows from
Lemma~\ref{lem-lim-D} that
\begin{eqnarray*}
\lim_{n\to \infty} [\vc{L}(n\tau+l)]_{i,j}
&=& [\vc{\psi}]_j \sum_{\nu=0}^{\tau-1}(\omega_{\tau}^{-\nu} )^{l-p(i)+p(j)}
\nonumber
\\
&=& \left\{
\begin{array}{ll}
\tau[\vc{\psi}]_j, & \mbox{if}~p(i) \equiv p(j)+l~(\Mod~\tau),
\\
0, & \mbox{otherwise}.
\end{array}
\right.
\end{eqnarray*}
This equation can be rewritten as
\begin{equation}
\lim_{n\to \infty} \vc{L}(n\tau+l)
= \tau \vc{E} \vc{H}_{l},
\label{eqn-lim-D(ntau+l)}
\end{equation}
where
\begin{equation}
\vc{E} = 
\bordermatrix{
              &           &         &           &        &             
\cr
\bbM^{(0)}        & \vc{e}^{(0)} & \vc{0}   & \cdots    & \vc{0} & \vc{0} 
\cr
\bbM^{(1)}        & \vc{0}   & \vc{e}^{(1)} & \cdots    & \vc{0} & \vc{0}
\cr
\ \vdots      & \vdots   & \vdots   & \ddots    & \vdots & \vdots
\cr
\bbM^{(\tau-2)} & \vc{0}   & \vc{0}   & \cdots    & \vc{e}^{(\tau-2)} & \vc{0}
\cr
\bbM^{(\tau-1)} & \vc{0}   & \vc{0}   &  \cdots   & \vc{0}   & \vc{e}^{(\tau-1)}
},
\label{defn-E}
\end{equation}
and
\[
\vc{H}_{l} = 
\bordermatrix{
              & \bbM^{(0)} & \bbM^{(1)} & \cdots & \bbM^{(\tau-l-1)} & \bbM^{(\tau-l)} 
&  \bbM^{(\tau-l+1)} & \cdots & \bbM^{(\tau-1)}         
\cr
& \vc{0} 
& \vc{0} 
& \cdots 
& \vc{0}  
& \vc{\psi}^{(\tau-l)} 
& \vc{0}
& \cdots
& \vc{0}
\cr
& \vc{0} 
& \vc{0} 
& \cdots 
& \vc{0} 
& \vc{0} 
& \vc{\psi}^{(\tau-l+1)} 
& \cdots
& \vc{0}
\cr
& \vdots
& \vdots
& \ddots
& \vdots 
& \vdots
& \vdots
& \ddots
& \vdots
\cr
& \vc{0} 
& \vc{0} 
& \cdots 
& \vc{0} 
& \vc{0} 
& \vc{0} 
& \cdots
& \vc{\psi}^{(\tau-1)} 
\cr
& \vc{\psi}^{(0)}
& \vc{0} 
& \cdots 
& \vc{0} 
& \vc{0} 
& \vc{0} 
& \cdots
& \vc{0}
\cr
& \vc{0}
& \vc{\psi}^{(1)}
& \cdots 
& \vc{0} 
& \vc{0} 
& \vc{0} 
& \cdots
& \vc{0}
\cr
& \vdots
& \vdots
& \ddots
& \vdots 
& \vdots
& \vdots
& \ddots
& \vdots
\cr
& \vc{0}
& \vc{0}
& \cdots 
& \vc{\psi}^{(\tau-l-1)} 
& \vc{0} 
& \vc{0} 
& \cdots
& \vc{0}
}.
\]

\medskip

\begin{rem}
Suppose the Markov chain $\{(X_n,S_n)\}$ is of M/G/1-type. It then
follows that $\vc{L}(n) = \vc{G}^n$ for $n=1,2,\dots$. Further it is
easy to see that $\vc{\psi}$ is a stationary probability vector of
$\vc{G}$ and therefore $[\vc{\psi}]_j = 0$ for all $j \in \bbM_{\rmT}$
(see Proposition~\ref{prop-structure-G}). We now define
$\vc{\psi}_{\bullet}^{(l)}$ ($l=0,1,\dots,\tau-1$) as a subvector of
$\vc{\psi}$ corresponding to $\bbM_{\bullet}^{(l)} :=
\{j\in\bbM_{\bullet} \cap \bbM^{(l)}\}$. As a result,
(\ref{eqn-lim-D(ntau+l)}) yields
\begin{eqnarray}
\lim_{n\to \infty} {1\over \tau}\vc{G}^{n\tau}
= \bordermatrix{
&   \bbM_{\bullet}^{(0)}  
& \bbM_{\bullet}^{(1)} 
& \cdots                 
& \bbM_{\bullet}^{(\tau-1)}             
& \bbM_{\rm T}
\cr
\bbM_{\bullet}^{(0)}        
& \vc{e}\vc{\psi}_{\bullet}^{(0)} 
& \vc{O}   
& \cdots    
& \vc{O}  
& \vc{O}
\cr
\bbM_{\bullet}^{(1)}        
& \vc{O}   
& \vc{e}\vc{\psi}_{\bullet}^{(1)}
& \cdots    
& \vc{O} 
& \vc{O}
\cr
\ \vdots      
& \vdots   
& \vdots   
& \ddots   
& \vdots 
& \vdots
\cr
\bbM_{\bullet}^{(\tau-1)} 
& \vc{O}   
& \vc{O}   
&  \cdots    
& \vc{e}\vc{\psi}_{\bullet}^{(\tau-1)} 
& \vc{O}
\cr
\bbM_{\rmT}^{(0)}        
& \vc{e}\vc{\psi}_{\bullet}^{(0)}  
& \vc{O}   
& \cdots     
& \vc{O}  
& \vc{O}
\cr
\bbM_{\rmT}^{(1)}        
& \vc{O}   
& \vc{e}\vc{\psi}_{\bullet}^{(1)}  
& \cdots    
& \vc{O} 
& \vc{O}
\cr
\ \vdots      
& \vdots   
& \vdots   
& \ddots    
& \vdots 
& \vdots
\cr
\bbM_{\rmT}^{(\tau-1)} 
& \vc{O}   
& \vc{O}   
&  \cdots     
& \vc{e}\vc{\psi}_{\bullet}^{(\tau-1)}   
& \vc{O}
},
\label{eqn-lim-G^{ntau}}
\end{eqnarray}
where $\bbM_{\rmT}^{(l)} = \bbM^{(l)} \setminus \bbM_{\bullet}^{(l)}$
($l=0,1,\dots,\tau-1$). Note here that
$\vc{\psi}_{\bullet}^{(l)}\vc{e}=1/\tau$ for all $l=0,1,\dots,\tau-1$
because $(1/ \tau)\vc{G}^{n\tau}\vc{e} = \vc{e}/\tau$ for all
$n=1,2,\dots$. As a result, the limit (\ref{eqn-lim-G^{ntau}}) is
consistent with the equation (14) in \cite{Masu11}, where
$\sum_{\nu=1}^{\tau}\vc{f}_{\nu}=\vc{e}$ and each element of
$\vc{f}_{\nu}$ ($\nu=1,2,\dots,\tau$) is equal to one or zero.
\end{rem}

\begin{lem}\label{lem-lim-sum-D(ntau+l)}
If Assumption~\ref{assu-sec3}~(I) holds, then
\begin{equation}
\lim_{n\to \infty} \sum_{l=0}^{\tau-1}\vc{L}(n\tau+l)
= \tau \vc{e} \vc{\psi}.
\label{eqn-lim-sum-D(ntau+l)}
\end{equation}
\end{lem}

\proof We obtain (\ref{eqn-lim-sum-D(ntau+l)}) by combining
(\ref{eqn-lim-D(ntau+l)}) and
\begin{equation}
\sum_{l=0}^{\tau-1}\vc{H}_l=\vc{e}\vc{\psi}.
\label{sum-H_l}
\end{equation}
\qed

\medskip

We now make the following assumption.
\begin{assumpt}\label{assu-sec3-A(k)e-B(k)e}
There exists some random variable $Y$ in $\bbZ_+$ with positive finite
mean such that
\begin{equation}
\lim_{k \to \infty} {\overline{\vc{A}}(k)\vc{e} \over \PP(Y > k)}
= {\vc{c}_{\!A} \over \EE[Y]},
\quad
\lim_{k \to \infty} {\overline{\vc{B}}(k)\vc{e} \over \PP(Y > k)}
= {\vc{c}_{\!B} \over \EE[Y]},
\label{eqn-lim-A(k)e}
\end{equation}
where $\vc{c}_{\!A}$ and $\vc{c}_{\!B}$ are nonnegative $M \times 1$
and $M_0 \times 1$ vectors, respectively, satisfying $\vc{c}_{\!A}
\neq \vc{0}$ or $\vc{c}_{\!B} \neq \vc{0}$.
\end{assumpt}

\begin{lem}\label{lem-AD-BD}
Suppose Assumptions~\ref{assu-sec3}~(I) and
\ref{assu-sec3-A(k)e-B(k)e} hold. If $Y_{\rm e}$ is long-tailed (i.e.,
$Y_{\rm e} \in \calL$; see Definition~\ref{defn-long-tailed}). We then
have
\begin{eqnarray}
\lim_{k \to \infty} \sum_{m = 1}^{\infty} 
{\overline{\vc{A}}(k+m) \vc{L}(m)
 \over \PP(Y_{\rm e} > k) }
&=& 
{\vc{c}_{\!A}\vc{\pi}(\vc{I} - \vc{R} ) (\vc{I} - \vc{\Phi}(0)) 
\over -\sigma },
\label{limit-AD}
\\
\lim_{k \to \infty} \sum_{m = 1}^{\infty} 
{\overline{\vc{B}}(k+m) \vc{L}(m)
 \over \PP(Y_{\rm e} > k) } 
&=& 
{\vc{c}_{\!B}\vc{\pi}(\vc{I} - \vc{R} ) (\vc{I} - \vc{\Phi}(0)) 
\over -\sigma }.
\label{limit-BD}
\end{eqnarray}
\end{lem}

\proof
See Appendix~\ref{proof-lem-AD=BD}. \qed

\begin{lem}\label{lem-R-R0}
Suppose Assumptions~\ref{assu-sec3}~(I) and
\ref{assu-sec3-A(k)e-B(k)e} hold. If $Y_{\rm e} \in \calL$, then
\begin{eqnarray}
\lim_{k \to \infty} {\overline{\vc{R}}(k) \over \PP(Y_{\rm e} > k) }
&=&
{\vc{c}_{\!A}\vc{\pi} ( \vc{I} - \vc{R} )  \over -\sigma
},
\label{eq-R-A-Y}
\\
\lim_{k \to \infty} {\overline{\vc{R}_0}(k) \over \PP(Y_{\rm e} > k) }
&=&
{\vc{c}_{\!B}\vc{\pi} ( \vc{I} - \vc{R} ) 
\over -\sigma
}.
\label{eq-R0-B-Y}
\end{eqnarray}
\end{lem}
\proof
It follows
from (\ref{eqn-R-A}) that
\begin{equation}
\overline{\vc{R}}(k) 
= \left[ 
\overline{\vc{A}}(k) + \sum_{m = 1}^{\infty} 
\overline{\vc{A}}(k+m) \vc{L}(m) \right](\vc{I} - \vc{\Phi}(0))^{-1}.
\label{eqn-overline{R}(k)}
\end{equation}
Note that Corollary 3.3 in \cite{Sigm99} and (\ref{eqn-lim-A(k)e}) yield
\[
\limsup_{k \to \infty} {\overline{\vc{A}}(k) \over \PP(Y_{\rm e} > k) }
\le \limsup_{k \to \infty} {\overline{\vc{A}}(k)\vc{e}\vc{e}^{\rm t} \over \PP(Y > k) }
\limsup_{k \to \infty} {\PP(Y > k) \over \PP(Y_{\rm e} > k) } 
= \vc{O}.
\]
Thus from (\ref{eqn-overline{R}(k)}), we have
\begin{equation}
\lim_{k \to \infty} {\overline{\vc{R}}(k) \over \PP(Y_{\rm e} > k) }
=
\lim_{k \to \infty} \sum_{m = 1}^{\infty} 
{\overline{\vc{A}}(k+m) \vc{L}(m)
 \over \PP(Y_{\rm e} > k) } ( \vc{I} - \vc{\Phi}(0) )^{-1} .
\label{eq-R-A-add1}
\end{equation}
Substituting (\ref{limit-AD}) into (\ref{eq-R-A-add1}), we obtain
(\ref{eq-R-A-Y}).  Similarly, we can prove (\ref{eq-R0-B-Y}). \qed

\medskip

The following theorem presents a subexponential asymptotic formula
for~$\{\overline{\vc{x}}(k)\}$.
\begin{thm}\label{thm-sub}
Suppose Assumptions~\ref{assu-sec3}~(I) and
\ref{assu-sec3-A(k)e-B(k)e} hold. If $Y_{\rm e} \in \calS$, then
\begin{equation}
\lim_{k\to \infty} {\overline{\vc{x}}(k) \over \PP(Y_{\rm e} > k) }
=
{\vc{x}(0)\vc{c}_{\!B} + \overline{\vc{x}}(0)\vc{c}_{\!A}
 \over
-\sigma
 }
 \cdot \vc{\pi}
\label{eq-bar-x-Y}
\end{equation}
\end{thm}

\proof 
It follows from (\ref{def-F(k)}) that
\begin{equation}
\sum_{k=0}^{\infty}\vc{F}(k) = (\vc{I} - \vc{R})^{-1}.
\label{eqn-sum-F(k)}
\end{equation}
Thus using (\ref{eqn-sum-F(k)}) and Lemma~6 in \cite{Jele98}, we have
\begin{eqnarray*}
\lim_{k \to \infty} {\overline{\vc{F}}(k) \over \PP(Y_{\rm e} > k)}
&=&
\lim_{k \to \infty} \sum_{n=0}^{\infty}{\overline{\vc{R}^{\ast n}}(k) 
\over \PP(Y_{\rm e} > k) } 
\\
&=&
 (\vc{I} - \vc{R} )^{-1} 
 \lim_{k \to \infty} {\overline{\vc{R}}(k) 
\over \PP(Y_{\rm e} > k) }(\vc{I} - \vc{R})^{-1}.
\end{eqnarray*}
Substituting (\ref{eq-R-A-Y}) into the above equation yields
\begin{equation}
\lim_{k \to \infty} {\overline{\vc{F}}(k) \over \PP(Y_{\rm e} > k)}
= 
 {(\vc{I} - \vc{R} )^{-1}
\vc{c}_{\!A}\vc{\pi} \over -\sigma
}.
\label{eq-gamma-Y}
\end{equation}
Finally, applying Proposition~A.3 in \cite{Masu11} to
(\ref{eq-xk-R0-Gamma}) and using (\ref{eq-R0-B-Y}) and
(\ref{eq-gamma-Y}) lead to
\[
\lim_{k \to \infty} {\overline{\vc{x}}(k) \over \PP(Y_{\rm e} > k) } 
=
{\vc{x}(0) \over -\sigma}
\left[ \vc{c}_{\!B} \vc{\pi}  +
\vc{R}_0(\vc{I} - \vc{R})^{-1} \vc{c}_{\!A} \vc{\pi}
\right],
\]
from which and (\ref{eq-bar-x_0-R}) we have (\ref{eq-bar-x-Y}). \qed

\begin{rem}
 Theorem~\ref{thm-sub} is a generalization of Theorem 1 in
 \cite{Kim11} to the GI/G/1-type Markov chain. In fact, the latter
 extends the corollary of Theorem 3.1 in \cite{Masu11} (Corollary~3.1
 therein) to the case where the $G$-matrix is periodic.
\end{rem}

\subsection{Case of strictly substochastic \vc{A}}\label{subsec-substochastic}

In this subsection, we make the following assumption in addition to
Assumption~\ref{assu-sec3}~(II):
\begin{assumpt}\label{assu-tail-A(k)-B(k)}
There exists some random variable $Y$ in $\bbZ_+$ such
that
\begin{equation}
\lim_{k \to \infty} {\overline{\vc{A}}(k) \over \PP(Y > k)}
= \vc{C}_{\!A},
\quad
\lim_{k \to \infty} {\overline{\vc{B}}(k) \over \PP(Y > k)}
= \vc{C}_{\!B},
\label{eq-assu-sub}
\end{equation}
where $\vc{C}_{\!A}$ and $\vc{C}_{\!B}$ are nonnegative $M \times M$
and $M_0 \times M$ matrices, respectively, satisfying $\vc{C}_{\!A}
\neq \vc{O}$ or $\vc{C}_{\!B} \neq \vc{O}$.
\end{assumpt}

\begin{lem}\label{lem-R-R0-sub}
Suppose Assumptions~\ref{assu-sec3}~(II) and \ref{assu-tail-A(k)-B(k)}
hold. If $Y \in \calL$, then
\begin{eqnarray}
\lim_{k \to \infty} {\overline{\vc{R}}(k) \over \PP(Y > k) }
&=&
\vc{C}_{\!A} 
\left( \vc{I} - \sum_{l=0}^{\infty} \vc{\Phi}(-l) \right)^{-1},
\label{eq-R-A-Y-sub}
\\
\lim_{k \to \infty} {\overline{\vc{R}_0}(k) \over \PP(Y > k) }
&=&
\vc{C}_{\!B}
\left( \vc{I} - \sum_{l=0}^{\infty} \vc{\Phi}(-l) \right)^{-1}.
\label{eq-R0-B-Y-sub}
\end{eqnarray}
\end{lem}
\proof
From (\ref{eqn-R-A}) and (\ref{eq-assu-sub}), we have
\begin{equation}
\lim_{k \to \infty}
{ \overline{\vc{R}}(k) \over \PP(Y > k) }
=
\left[
\vc{C}_A
 + \lim_{k \to \infty}\sum_{m=1}^{\infty} { \overline{\vc{A}}(k+m) \over \PP(Y > k)} 
\vc{L}(m)\right]
(\vc{I} - \vc{\Phi}(0))^{-1}. 
\label{eqn-R-sub}
\end{equation}
Note here that under Assumption~\ref{assu-sec3}~(II), $\rsp(\vc{G}) < 1$
(see Proposition~\ref{prop-Phi}) and thus (\ref{defn-hat{L}(z)})
yields
\begin{equation}
\sum_{m=1}^{\infty}\vc{L}(m) = (\vc{I} - \vc{G} )^{-1} \vc{G} < \infty,
\label{sum-L(m)}
\end{equation}
from which and (\ref{eq-assu-sub}) it follows that for $k=0,1,\dots$,
\[
\sum_{m=1}^{\infty}
{ \overline{\vc{A}}(k+m) \over \PP(Y > k)} \vc{L}(m)
\le 
\sup_{k\in\bbZ_+}{ \overline{\vc{A}}(k) \over \PP(Y > k)}
\sum_{m=1}^{\infty} \vc{L}(m) < \infty.
\]
Therefore applying the dominated convergence theorem to
(\ref{eqn-R-sub}) and using (\ref{eq-assu-sub}) and $Y \in \calL$, we
obtain
\begin{eqnarray}
\lefteqn{
\lim_{k \to \infty}
{ \overline{\vc{R}}(k) \over \PP(Y > k) }
}
~~&&
\nonumber
\\
&=&
\left[
\vc{C}_A
 + \sum_{m=1}^{\infty} 
\lim_{k \to \infty}{ \overline{\vc{A}}(k+m) \over \PP(Y > k+m)} 
{ \PP(Y > k+m) \over \PP(Y > k)} 
\vc{L}(m)\right]
(\vc{I} - \vc{\Phi}(0))^{-1}
\nonumber
\\
&=& \vc{C}_{\!A}\left[\vc{I} + (\vc{I} - \vc{G})^{-1}\vc{G} \right]
(\vc{I} - \vc{\Phi}(0))^{-1}
\nonumber
\\
&=& \vc{C}_{\!A}
(\vc{I} - \vc{G})^{-1}
(\vc{I} - \vc{\Phi}(0))^{-1}.
\label{eqn-R-sub-2}
\end{eqnarray}
From (\ref{def-G(k)}), we have
\begin{eqnarray}
(\vc{I} - \vc{G} )^{-1} 
&=&
\left[\vc{I} - (\vc{I} - \vc{\Phi}(0))^{-1}
\sum_{l=1}^{\infty} \vc{\Phi}(-l) \right]^{-1}
\nonumber
\\
&=&
\left( \vc{I} - \sum_{l=0}^{\infty} \vc{\Phi}(-l) \right)^{-1}
(\vc{I} - \vc{\Phi}(0)).
\label{add-eqn-03}
\end{eqnarray}
Finally, substituting (\ref{add-eqn-03}) into (\ref{eqn-R-sub-2})
yields (\ref{eq-R-A-Y-sub}). Equation (\ref{eq-R0-B-Y-sub}) can be
proved in the same way.  \qed

\setcounter{thmthm}{1}
\begin{thmthm}\label{thm-substoch}
Suppose Assumptions~\ref{assu-sec3}~(II) and \ref{assu-tail-A(k)-B(k)}
hold. If $Y \in \calS$, then
\begin{equation}
\lim_{k\to \infty} {\overline{\vc{x}}(k) \over \PP(Y > k) }
=
[
\vc{x}(0) \vc{C}_{\!B} + \overline{\vc{x}}(0) \vc{C}_{\!A}
] (\vc{I}- \vc{A})^{-1} > \vc{0}.
\label{eq-bar-x-Y-01}
\end{equation}
\end{thmthm}

\proof Applying Proposition~A.3 in \cite{Masu11} to
(\ref{eq-xk-R0-Gamma}) and using (\ref{eqn-sum-F(k)}) and
(\ref{eq-R0-B-Y-sub}), we have
\begin{eqnarray}
\lim_{k\to \infty} {\overline{\vc{x}}(k) \over \PP(Y > k) }
&=& \vc{x}(0) \vc{C}_{\!B}
\left( \vc{I} - \sum_{l=0}^{\infty} \vc{\Phi}(-l) \right)^{-1}
(\vc{I}- \vc{R})^{-1}
\nonumber
\\
&& {} + \vc{x}(0)\vc{R}_0 
\lim_{k\to \infty} {\overline{\vc{F}}(k) \over \PP(Y > k)},
\label{add-eqn-10}
\end{eqnarray}
where $\vc{F}(k)$ is given in (\ref{def-F(k)}). Further it follows
from Lemma~6 in \cite{Jele98} and
(\ref{eq-R-A-Y-sub}) that
\[
\lim_{k \to \infty} {\overline{\vc{F}}(k) \over \PP(Y > k)}
= 
(\vc{I} - \vc{R})^{-1}\vc{C}_{\!A}
\left( \vc{I} - \sum_{l=0}^{\infty} \vc{\Phi}(-l) \right)^{-1} 
(\vc{I} - \vc{R})^{-1}.
\]
Substituting the above equation into (\ref{add-eqn-10}) and using
(\ref{eq-bar-x_0-R}), we have
\begin{equation}
\lim_{k\to \infty} {\overline{\vc{x}}(k) \over \PP(Y > k) }
=
[
\vc{x}(0) \vc{C}_{\!B} + \overline{\vc{x}}(0) \vc{C}_{\!A}
]
 \left( \vc{I} -
 \sum_{l=0}^{\infty} \vc{\Phi}(-l) \right)^{-1} (\vc{I}- \vc{R})^{-1}.
\label{eq-bar-x-Y-3}
\end{equation}
Note here that (\ref{add-eqn-03}) yields
\begin{eqnarray}
\lefteqn{
\left( \vc{I} - \sum_{l=0}^{\infty} \vc{\Phi}(-l) \right)^{-1} 
(\vc{I} - \vc{R})^{-1}
}
\quad && \nonumber
\\
&=& \left( \vc{I} - \vc{G} \right)^{-1} 
(\vc{I} - \vc{\Phi}(0))^{-1}(\vc{I} - \vc{R})^{-1}
= (\vc{I} - \vc{A})^{-1},
\label{add-eqn-08}
\end{eqnarray}
where the second equality follows from Proposition~\ref{prop-RG}. As a
result, we obtain (\ref{eq-bar-x-Y-01}) by combining
(\ref{eq-bar-x-Y-3}) with (\ref{add-eqn-08}).

It is easy to show that the right hand side of (\ref{eq-bar-x-Y-01})
is positive. Indeed, $(\vc{I} - \vc{A})^{-1} > \vc{O}$ due to the
irreducibility of $\vc{A}$. In addition, $\vc{x}(0) \vc{C}_{\!B} +
\overline{\vc{x}}(0) \vc{C}_{\!A} \ge \vc{0},\neq \vc{0}$ because
$\vc{x}(0) > \vc{0}$ and $\overline{\vc{x}}(0) > \vc{0}$; and
$\vc{C}_{\!A} \neq \vc{O}$ or $\vc{C}_{\!B} \neq \vc{O}$. Therefore,
$(\vc{x}(0) \vc{C}_{\!B} + \overline{\vc{x}}(0) \vc{C}_{\!A})(\vc{I} -
\vc{A})^{-1} > \vc{0}$.  \qed

\section{Locally Subexponential Asymptotics}\label{sec-loc-subexp}

This section considers the locally subexponential asymptotics of the
stationary distribution.

\subsection{Case of stochastic \vc{A}}

In this subsection, we proceed under Assumption~\ref{assu-sec3}~(I)
and the following assumption:
\begin{assumpt}\label{assu-sec4-A(k)E-B(k)E}
There exists some random variable $Y$ in $\bbZ_+$ with positive finite
mean such that
\begin{equation}
\lim_{k \to \infty} {\vc{A}(k)\vc{E} \over \PP(Y = k)}
= {\vc{C}{}_{\!A}^{\rmE} \over \EE [Y]},
\quad
\lim_{k \to \infty} {\vc{B}(k)\vc{E} \over \PP(Y = k)}
= {\vc{C}{}_{\!B}^{\rmE} \over \EE [Y]},
\label{sec4-eq-assu-AE-BE}
\end{equation}
where $\vc{E}$ is given in (\ref{defn-E}), and where
$\vc{C}_{\!A}^{\rmE}$ and $\vc{C}_{\!B}^{\rmE}$ are nonnegative $M
\times \tau$ and $M_0 \times \tau$ matrices, respectively, satisfying
$\vc{C}_{\!A}^{\rmE} \neq \vc{O}$ or $\vc{C}_{\!B}^{\rmE} \neq
\vc{O}$.
\end{assumpt}

\begin{lem}\label{sec4-lem-AD-BD}
Suppose Assumptions~\ref{assu-sec3}~(I) and
\ref{assu-sec4-A(k)E-B(k)E} hold. Further, suppose either of the
following is satisfied: $Y$ is {\it locally long-tailed with span one}
(i.e., $Y \in \calL_{\rm loc}(1)$; see Definition~\ref{defn-loc-L});
or $Y \in \calL$ and $\{\PP(Y = k)\}$ is eventually
nonincreasing. Then
\begin{eqnarray}
\lim_{k \to \infty} \sum_{m = 1}^{\infty} 
{\vc{A}(k+m) \vc{L}(m)
 \over \PP(Y_{\rme} = k) }
&=& \vc{C}{}_{\!A}^{\rmE} \vc{e}
{\vc{\pi}(\vc{I} - \vc{R} ) (\vc{I} - \vc{\Phi}(0)) 
\over -\sigma },
\label{sec4-limit-AD}
\\
\lim_{k \to \infty} \sum_{m = 1}^{\infty} 
{\vc{B}(k+m) \vc{L}(m)
 \over \PP(Y_{\rme} = k) } 
&=& 
\vc{C}{}_{\!B}^{\rmE} \vc{e}
{\vc{\pi}(\vc{I} - \vc{R} ) (\vc{I} - \vc{\Phi}(0)) 
\over -\sigma }.
\label{sec4-limit-BD}
\end{eqnarray}
\end{lem}
\proof See Appendix~\ref{proof-sec4-lem-AD-BD}. \qed

\begin{rem}
Lemma~\ref{sec4-lem-AD-BD} is proved by using
Proposition~\ref{add-appen-prop2}, which requires either that $Y \in
\calL_{\rm loc}(1)$ or that $Y \in \calL$ and $\{\PP(Y = k)\}$ is
eventually nonincreasing.
\end{rem}

\begin{lem}\label{sec4-lem-R-R0}
Under the same assumptions as in Lemma~\ref{sec4-lem-AD-BD}, 
\begin{eqnarray}
\lim_{k \to \infty} {\vc{R}(k) \over\PP(Y_{\rme} = k) }
&=&
\vc{C}{}_{\!A}^{\rmE}\vc{e} {\vc{\pi} ( \vc{I} - \vc{R} )  \over -\sigma
},
\label{sec4-eq-R-A-Y}
\\
\lim_{k \to \infty} {\vc{R}_0(k) \over \PP(Y_{\rme} = k) }
&=&
\vc{C}{}_{\!B}^{\rmE}\vc{e} {\vc{\pi} ( \vc{I} - \vc{R} ) 
\over -\sigma
}.
\label{sec4-eq-R0-B-Y}
\end{eqnarray}
\end{lem}
\proof It follows from $\vc{E}\vc{e}=\vc{e}$,
(\ref{sec4-eq-assu-AE-BE}) and $Y \in \calL$ that
\[
\lim_{k \to \infty} {\vc{A}(k) \over \PP(Y_{\rme} = k) }
\le \EE[Y]\lim_{k \to \infty} 
{\vc{A}(k)\vc{E}\vc{e}\vc{e}^{\rmt} \over \PP(Y > k) }
{\PP(Y=k) \over \PP(Y > k) }
= \vc{O}. 
\]
Thus from (\ref{eqn-R-A}), we have
\begin{equation}
\lim_{k \to \infty} {\vc{R}(k) \over \PP(Y_{\rme} = k) }
=
\lim_{k \to \infty} \sum_{m = 1}^{\infty} 
{\vc{A}(k+m) \vc{L}(m) \over \PP(Y_{\rme} = k) } 
( \vc{I} - \vc{\Phi}(0) )^{-1} .
\label{sec4-eq-R-A-add1}
\end{equation}
Substituting (\ref{sec4-limit-AD}) into (\ref{sec4-eq-R-A-add1})
yields (\ref{sec4-eq-R-A-Y}). Similarly, we can readily show
(\ref{sec4-eq-R0-B-Y}). \qed

\medskip

We now obtain a locally subexponential asymptotic formula for
$\{\vc{x}(k)\}$.
\begin{thm}\label{sec4-thm-sub}
Suppose Assumptions~\ref{assu-sec3}~(I) and
\ref{assu-sec4-A(k)E-B(k)E} hold. Further, suppose (i) $Y_{\rme}$ is
    {\it locally subexponential with span one} (i.e., $Y_{\rme} \in
    \calS_{\rm loc}(1)$; see Definition~\ref{defn-loc-S}); and (ii) $Y
    \in \calL_{\rm loc}(1)$ or $\{\PP(Y = k)\}$ is eventually
    nonincreasing. Then
\begin{equation}
\lim_{k\to \infty} {\vc{x}(k) \over \PP(Y_{\rme} = k) }
=
{\vc{x}(0)
 \vc{C}{}_{\!B}^{\rmE} \vc{e} + \overline{\vc{x}}(0)\vc{C}{}_{\!A}^{\rmE} \vc{e}
 \over -\sigma
 }
 \cdot \vc{\pi}.
\label{sec4-eq-bar-x-Y}
\end{equation}
\end{thm}

\begin{rem}
According to Definition~\ref{defn-loc-S} and
Proposition~\ref{prop-S^{ast}-S_{loc}}, $Y_{\rme} \in \calS_{\rm
  loc}(1)$ is equivalent to $Y \in \calS^{\ast}$. Thus since
$\calS^{\ast} \subset \calS \subset \calL$, the assumptions of
Theorem~\ref{sec4-thm-sub} are sufficient for those of
Lemma~\ref{sec4-lem-AD-BD}.

\end{rem}
 
\noindent
{\it Proof of Theorem~\ref{sec4-thm-sub}}.~
Proposition~\ref{appen-prop-convo2} yields
\begin{eqnarray*}
\lim_{k \to \infty} {\vc{F}(k) \over \PP(Y_{\rme} = k)}
&=&
\lim_{k \to \infty} 
\sum_{n=0}^{\infty}{\vc{R}^{\ast n}(k) \over \PP(Y_{\rm e} = k) } 
\\
&=& 
 (\vc{I} - \vc{R} )^{-1} 
 \lim_{k \to \infty} {\vc{R}(k) 
\over \PP(Y_{\rm e} = k) }(\vc{I} - \vc{R})^{-1},
\end{eqnarray*}
from which and (\ref{sec4-eq-R-A-Y}) it follows that
\begin{equation}
\lim_{k \to \infty} {\vc{F}(k) \over \PP(Y_{\rme} = k)}
= 
 {(\vc{I} - \vc{R} )^{-1}
\vc{C}{}_{\!A}^{\rmE}\vc{e} \vc{\pi} \over -\sigma
}.
\label{sec4-eq-gamma-Y}
\end{equation}
Further applying Proposition~\ref{appen-prop-convo1} to
(\ref{eqn-x(k)-convolution}) and using (\ref{sec4-eq-R0-B-Y}) and
(\ref{sec4-eq-gamma-Y}), we obtain
\[
\lim_{k \to \infty} {\vc{x}(k) \over \PP(Y_{\rm e} = k) } 
= 
{\vc{x}(0) \over -\sigma}
\left[ \vc{C}{}_{\!B}^{\rmE} \vc{e} \vc{\pi}  +
\vc{R}_0(\vc{I} - \vc{R})^{-1} \vc{C}{}_{\!A}^{\rmE} \vc{e} \vc{\pi}
\right].
\]
Substituting (\ref{eq-bar-x_0-R}) into the above equation yields
(\ref{sec4-eq-bar-x-Y}). \qed

\medskip

We present another asymptotic formula.
\begin{assumpt}\label{sec4-assu-A(k)e-B(k)e}
There exists some random variable
  $Y$ in $\bbZ_+$ with positive finite mean such that
\[
\lim_{k \to \infty} {\vc{A}(k)\vc{e} \over \PP(Y = k)}
={ \vc{c}_{\!A} \over \EE [Y]},
\quad
\lim_{k \to \infty} {\vc{B}(k)\vc{e} \over \PP(Y = k)}
= { \vc{c}_{\!B} \over \EE [Y]},
\]
where $\vc{c}_{\!A}$ and $\vc{c}_{\!B}$ are nonnegative $M \times 1$
and $M_0 \times 1$ vectors, respectively, satisfying $\vc{c}_{\!A}
\neq \vc{0}$ or $\vc{c}_{\!B} \neq \vc{0}$.
\end{assumpt}
\begin{thm}\label{sec4-thm-sub2}
Suppose Assumptions~\ref{assu-sec3}~(I) and
\ref{sec4-assu-A(k)e-B(k)e} hold.  Further, suppose (i) $Y_{\rme} \in
\calS_{\rm loc}(1)$; (ii) $Y \in \calL_{\rm loc}(1)$ or $\{\PP(Y =
k)\}$ is eventually nonincreasing; and (iii)
$\{\vc{A}(k);k\in\bbZ_+\}$ and $\{\vc{B}(k);k\in\bbN\}$ are eventually
nonincreasing. Then
\[
\lim_{k\to \infty} {\vc{x}(k) \over \PP(Y_{\rme} = k) }
=
{\vc{x}(0)\vc{c}_{\!B} + \overline{\vc{x}}(0)\vc{c}_{\!A}
 \over -\sigma
 }
 \cdot \vc{\pi}.
\]
\end{thm}

\proof This theorem can be proved in a very similar way to
Theorem~\ref{thm-sub}.  For doing this, we require an
additional condition that $\{\vc{A}(k);k\in\bbZ_+\}$ and
$\{\vc{B}(k);k\in\bbN\}$ are eventually nonincreasing, i.e., there
exists some $k_{\ast} \in \bbN$ such that $\vc{A}(k) \ge \vc{A}(k+1)$
and $\vc{B}(k) \ge \vc{B}(k+1)$ for all $k \ge k_{\ast}$.  The details
are omitted.  \qed

\begin{rem}
Since $\vc{E}\vc{e} = \vc{e}$, Assumption~\ref{sec4-assu-A(k)e-B(k)e}
is sufficient for Assumption~\ref{assu-sec4-A(k)E-B(k)E}. Thus
Theorem~\ref{sec4-thm-sub2} is not a collorary of
Theorem~\ref{sec4-thm-sub}.
\end{rem}

\subsection{Case of strictly substochastic \vc{A}}

In addition to Assumption~\ref{assu-sec3}~(II), we assume the
following:
\begin{assumpt}\label{sec4-assu-sub}
There exists some random variable $Y$ in $\bbZ_+$ such that
\begin{equation}
\lim_{k \to \infty} {\vc{A}(k) \over \PP(Y = k)}
=\vc{C}_{\!A},
\quad
\lim_{k \to \infty} {\vc{B}(k) \over \PP(Y = k)}
=\vc{C}_{\!B},
\label{sec4-eq-assu-sub}
\end{equation}
where $\vc{C}_{\!A}$ and $\vc{C}_{\!B}$ are nonnegative $M \times M$
and $M_0 \times M$ matrices, respectively, satisfying $\vc{C}_{\!A}
\neq \vc{O}$ or $\vc{C}_{\!B} \neq \vc{O}$.
\end{assumpt}

\begin{lem}\label{sec4-lem-R-R0-sub}
Suppose Assumptions~\ref{assu-sec3}~(II) and \ref{sec4-assu-sub} hold.
If $Y \in \calL_{\rm loc}(1)$; and $r_{A_-} > 1$ or $\{\PP(Y = k)\}$
is eventually nonincreasing, then
\begin{eqnarray}
\lim_{k \to \infty} {\vc{R}(k) \over \PP(Y = k) }
&=&
\vc{C}_{\!A} 
\left( \vc{I} - \sum_{l=0}^{\infty} \vc{\Phi}(-l) \right)^{-1},
\label{sec4-eq-R-A-Y-sub}
\\
\lim_{k \to \infty} {\vc{R}_0(k) \over \PP(Y = k) }
&=&
\vc{C}_{\!B}
\left( \vc{I} - \sum_{l=0}^{\infty} \vc{\Phi}(-l) \right)^{-1}.
\label{sec4-eq-R0-B-Y-sub}
\end{eqnarray}
\end{lem}

\proof From (\ref{eqn-R-A}) and (\ref{sec4-eq-assu-sub}), we have
\begin{equation}
\lim_{k \to \infty}
{ \vc{R}(k) \over \PP(Y = k) }
=
\left[
\vc{C}_{\!A} 
 + \lim_{k \to \infty}\sum_{m=1}^{\infty} { \vc{A}(k+m) \over \PP(Y = k)} 
\vc{L}(m)\right]
(\vc{I} - \vc{\Phi}(0))^{-1}. 
\label{eqn-R-sub'}
\end{equation}
To apply the dominated convergence theorem to (\ref{eqn-R-sub'}), we
show that for all sufficiently large $k$,
\[
\sum_{m=1}^{\infty}
{ \vc{A}(k+m) \over \PP(Y = k)} \vc{L}(m) < \infty.
\]
Suppose $\{\PP(Y = k)\}$ is eventually nonincreasing. We then have for
all sufficiently large $k$,
\[
\sum_{m=1}^{\infty}
{ \vc{A}(k+m) \over \PP(Y = k)} \vc{L}(m)
\le 
\sup_{m'\in\bbN}{ \vc{A}(k+m') \over \PP(Y = k+m')}
\sum_{m=1}^{\infty} \vc{L}(m)<\infty,
\]
where the last inequality is due to (\ref{sum-L(m)}) and
(\ref{sec4-eq-assu-sub}).  On the other hand, suppose $r_{A_-} >
1$. It then follows from Proposition~\ref{prop-convergence-radii} that
$\{\vc{G}(k)\}$ is light-tailed, i.e.,
\begin{equation}
\sum_{k=1}^{\infty}r^k\vc{G}(k) < \infty\quad \mbox{for all}~1 < r < r_{A_-}.
\label{add-eqn-33}
\end{equation}
Note here that $\widehat{\vc{G}}(1/z) =
\sum_{k=1}^{\infty}z^k\vc{G}(k)$ and $\rsp(\widehat{\vc{G}}(1)) < 1$
(see Proposition~\ref{prop-Phi}). Thus according to Theorem~8.1.18 in
\cite{Horn90},
\begin{equation}
\rsp(\widehat{\vc{G}}(1/z)) = 1\quad \mbox{only if}~ 1 < z \le r_{A_-}.
\label{add-eqn-34}
\end{equation}
The equations (\ref{defn-hat{L}(z)}), (\ref{add-eqn-33}) and
(\ref{add-eqn-34}) imply that there exists some $r > 1$ such that
\[
\sum_{m=1}^{\infty}r^m \vc{L}(m) < \infty.
\]
Further it follows from Assumption~\ref{sec4-assu-sub} and $Y \in
\calL_{\rm loc}(1)$ that for any $\varepsilon > 0$ there exists some
$k_0 \in \bbZ_+$ such that for all $k \ge k_0$,
\[
{ \vc{A}(k+m) \over \PP(Y = k)} 
\le (\vc{C}_{\!A} + \varepsilon \vc{e} \vc{e}^{\rmt}) { \PP(Y = k+m) \over \PP(Y = k)}
\le (1 + \varepsilon)^{m} (\vc{C}_{\!A} + \varepsilon \vc{e} \vc{e}^{\rmt}),
\quad m \in \bbZ_+.
\]
Therefore, for $0 < \varepsilon \le r - 1$ and $k \ge k_0$,
\[
\sum_{m=1}^{\infty}
{ \vc{A}(k+m) \over \PP(Y = k)} \vc{L}(m)
\le 
(\vc{C}_{\!A} + \varepsilon \vc{e} \vc{e}^{\rmt})
\sum_{m=1}^{\infty} (1 + \varepsilon)^{m}\vc{L}(m)<\infty.
\]

As a result, applying the dominated convergence theorem to
(\ref{eqn-R-sub'}) and following the proof of
Lemma~\ref{lem-R-R0-sub}, we can prove
(\ref{sec4-eq-R-A-Y-sub}). Equation (\ref{sec4-eq-R0-B-Y-sub}) can be
proved in the same way.  \qed

\medskip

Using Lemma~\ref{sec4-lem-R-R0-sub}, we can readily prove the
following theorem. The proof is very similar to that of
Theorem~\ref{thm-substoch} and thus is omitted.
\begin{thm}\label{thm-loc-substoch}
Suppose Assumptions~\ref{assu-sec3}~(II) and \ref{sec4-assu-sub} hold.
If $Y \in \calS_{\rm loc}(1)$; and $r_{A_-} > 1$ or $\{\PP(Y = k)\}$
is eventually nonincreasing, then
\begin{eqnarray*}
\lim_{k\to \infty} {\vc{x}(k) \over \PP(Y = k) }
&=&
[
\vc{x}(0) \vc{C}_{\!B} + \overline{\vc{x}}(0) \vc{C}_{\!A}
]
 (\vc{I}- \vc{A})^{-1} > \vc{0}.
\end{eqnarray*}
\end{thm}

\section{Discussion on Assumptions}

This section discusses the assumptions of the theorems presented in
Sections~\ref{sec-subexp} and \ref{sec-loc-subexp}.

We first consider the case of stochastic $\vc{A}$, for which
Theorems~\ref{thm-sub}, \ref{sec4-thm-sub} and \ref{sec4-thm-sub2} are
shown. The assumptions of these theorems are summarized in
Table~\ref{tab:stochastic}, where ``{\it eventually nonincreasing}" is
abbreviated as ``ENI".
\begin{table}[htb]
\begin{center}
\caption{The assumptions of the theorems in case of stochastic $\vc{A}$}
\begin{tabular}{|r||r||r|} \hline
  Theorem \ref{thm-sub} 
& Theorem \ref{sec4-thm-sub}
& Theorem \ref{sec4-thm-sub2}
\\ 
\hline \hline 
  Assumption \ref{assu-sec3} (I) 
& Assumption \ref{assu-sec3} (I) 
& Assumption \ref{assu-sec3} (I) 
\\ 
\hline 
  Assumption \ref{assu-sec3-A(k)e-B(k)e} 
& Assumption \ref{assu-sec4-A(k)E-B(k)E} 
& Assumption \ref{sec4-assu-A(k)e-B(k)e} 
\\ 
\hline

  $Y_{\rm e} \in \calS$ 
& $Y_{\rm e} \in \calS_{\rm loc}(1)$ 
& $Y_{\rm e} \in \calS_{\rm loc}(1)$ 
\\ 
\hline 
& $Y \in \calL_{\rm loc}(1)$ or 
& $Y \in \calL_{\rm loc}(1)$ or 
\\ 
&
  $\{\PP(Y=k)\}$ is ENI 
& $\{\PP(Y=k)\}$ is ENI 
\\ 
\hline & & $\{\vc{A}(k)\}$ and $\{\vc{B}(k)\}$ are ENI 
\\ \hline
\end{tabular}
\label{tab:stochastic}
\end{center}
\end{table}
Note here that Assumption~\ref{assu-sec4-A(k)E-B(k)E} implies
Assumption~\ref{assu-sec3-A(k)e-B(k)e} due to $\vc{E}\vc{e} =
\vc{e}$. Recall also that if $Y_{\rm e} \in \calS_{\rm loc}(1)$, then
$Y_{\rm e} \in \calS$ (see Remark~\ref{rem-local-subexp}). Thus the
assumptions of Theorem~\ref{sec4-thm-sub} are more restrictive than
those of Theorem~\ref{thm-sub}.  Similarly, we can readily confirm
that the assumptions of Theorem~\ref{sec4-thm-sub2} imply those of
Theorem~\ref{thm-sub}.  It should be noted that
Theorem~\ref{sec4-thm-sub2} is not a corollary of
Theorem~\ref{sec4-thm-sub} because
Assumption~\ref{sec4-assu-A(k)e-B(k)e} is weaker than
Assumption~\ref{assu-sec4-A(k)E-B(k)E}.

Next we consider the case of substochastic $\vc{A}$, for which
Theorems~\ref{thm-substoch} and \ref{thm-loc-substoch} are shown.  It
is easy to see that Assumption~\ref{sec4-assu-sub} implies
Assumption~\ref{assu-tail-A(k)-B(k)}. Further if $Y \in \calS_{\rm
  loc}(1)$, then $Y \in \calS$ (see
Remark~\ref{rem-local-subexp}). Therefore the assumptions of
Theorem~\ref{thm-loc-substoch} are more restrictive than those of
Theorem~\ref{thm-substoch} (see Table~\ref{tab:substochastic}).
\begin{table}[htb]
\begin{center}
\caption{The assumptions of the theorems in case of strictly substochastic $\vc{A}$}
\begin{tabular}{|r||r|} \hline
Theorem \ref{thm-substoch} & Theorem \ref{thm-loc-substoch}
\\ 
\hline \hline
  Assumption \ref{assu-sec3} (II)  
& Assumption \ref{assu-sec3} (II) 
\\ 
\hline
  Assumption \ref{assu-tail-A(k)-B(k)}  
& Assumption \ref{sec4-assu-sub} 
\\ 
\hline
  $Y \in \calS$ 
& $Y \in \calS_{\rm loc}(1)$
\\ 
\hline
 & $r_{A_-} > 1$ or
\\ 
 & $\{\PP(Y=k)\}$ is ENI
\\ 
\hline
\end{tabular}
\label{tab:substochastic}
\end{center}
\end{table}

\appendix

\section{Subexponential Distributions}

This section provides a brief overview of two classes of
subexponential distributions on $\bbZ_+$. One is the class of
``ordinal" subexponential distributions introduced by Chistyakov
\cite{Chis64}, and the other one is the class of ``locally"
subexponential distributions introduced by Chover et al.~\cite{Chov73}
and generalized by Asmussen et al.~\cite{Asmu03}.

In what follows, let $U$ denote a random variable in $\bbZ_+$ and
$U_j$ ($j\in\bbZ_+$) denote independent copies of $U$. Let $U_{\rm e}$
denote the discrete equilibrium random variable of $U$, distributed
with $\PP(U_{\rm e} = k) = \PP(U > k)/\EE[U]$ ($k\in\bbZ_+$). Further,
for any $h \in \bbN \cup\{\infty\}$, let $\Delta_h = (0,h]$ and
  $k+\Delta_h = \{x \ge0; k < x \le k+h\}$ for $k \in \bbZ_+$.

\subsection{Ordinal subexponential class}

We begin with the definition of the long-tailed class, which covers
the subexponential class.
\begin{app-defn}[\cite{Asmu03-book,Embr97,Sigm99}]\label{defn-long-tailed}
A random variable $U$ in $\bbZ_+$ and its distribution are said to be
long-tailed if $\PP(U>k) > 0$ for all $k\in\bbZ_+$ and
\[
\lim_{k\to\infty}{\PP(U > k+1) \over \PP(U>k)} = 1.
\]
The class of long-tailed distributions is denoted by $\calL$.
\end{app-defn}

The following result is used to derive some of the asymptotic results
presented in this paper.
\begin{app-prop-except}[Proposition~A.1 in \cite{Masu11}]\label{add-appen-prop}
If $U_{\rm e} \in \calL$, then for any $h \in \bbN$, $l_0\in \bbZ_+$ and
$\nu=0,1,\dots,h-1$,
\[
{1 \over \EE[U]}
\lim_{k\to\infty}
{\sum_{l=l_0}^{\infty}\PP(U > k + lh + \nu) \over \PP(U_{\rm e} > k)}
= {1 \over h}.
\]
\end{app-prop-except}

We now introduce the definition of the subexponential class.
\begin{app-defn}[\cite{Chis64,Embr97,Sigm99}]\label{defn-subexp}
A random variable $U$ and its distribution are said to be
subexponential if $\PP(U>k) > 0$ for all $k\in\bbZ_+$ and
\[
\lim_{k\to\infty}{\PP(U_1 + U_2 > k) \over \PP(U>k)} = 2.
\]
The class of subexponential
distributions is denoted by $\calS$.
\end{app-defn}

\begin{app-rem}
$\calS \subset \calL$ (see, e.g., \cite{Sigm99}), and there exists an
  example of not subexponential but long-tailed distributions (see
  \cite{Pitm80}).
\end{app-rem}

The following is a discrete analog of class $\calS^{\ast}$ introduced
by Kl\"{u}ppelberg~\cite{Klup88}.

\begin{app-defn}\label{def-S^*}
A random variable $U$ and its distribution belong to class
$\calS^{\ast}$ if $\PP(U>k) > 0$ for all $k\in\bbZ_+$ and
\begin{equation}
\lim_{k\to\infty}
\sum_{l=0}^k {\PP(U > k-l)\PP(U > l) \over \PP(U > k)}
= 2\EE[U] < \infty.
\label{def-S*}
\end{equation}
\end{app-defn}

\begin{app-rem}
 If $U \in \calS^{\ast}$, then $U \in \calS$ and $U_{\rm e} \in \calS$
 (see Proposition~A.2 in \cite{Masu11}).
\end{app-rem}

\subsection{Locally subexponential class}

We first introduce the locally long-tailed class, which is required by
the definition of the locally subexponential class.
\begin{app-defn}[Definition~1 in \cite{Asmu03}]\label{defn-loc-L}
A random variable $U$ and its distribution $F$ are called {\it locally
  long-tailed} with span $h \in \bbN \cup \{\infty\}$ if $\PP(U \in
k+\Delta_h) > 0$ for all sufficiently large $k$ and
\[
\lim_{k\to\infty}
{\PP(U \in k+1+\Delta_h) \over \PP(U \in k+\Delta_h)} 
= 1.
\]
\end{app-defn}

We denote by $\calL_{\rm loc}(h)$ the class of locally long-tailed
distributions with span $h$ hereafter.

\begin{app-rem}\label{rem-L_loc}
By definition, $\calL_{\rm loc}(\infty) = \calL$. Further, if $U \in
\calL_{\rm loc}(1)$, then $U \in \calL_{\rm loc}(n)$ for all
$n=2,3,\dots$ and $U \in \calL$.
\end{app-rem}

The following proposition is a locally asymptotic version of
Proposition~\ref{add-appen-prop}.
\begin{app-prop}\label{add-appen-prop2}
Suppose (i) $U \in \calL_{\rm loc}(1)$; or (ii) $U \in \calL$ and
$\{\PP(U=k)\}$ is eventually nonincreasing. Then for any $h \in \bbN$,
$l_0\in \bbZ_+$ and $\nu=0,1,\dots,h-1$,
\begin{equation}
\lim_{k\to\infty}
{\sum_{l=l_0}^{\infty}\PP(U = k + lh + \nu) \over \PP(U > k)}
= {1 \over h}.
\label{add-eqn-24}
\end{equation}
\end{app-prop}

\proof See Appendix~\ref{proof-add-appen-prop2}. \qed

\begin{app-defn}[Definition~2 in \cite{Asmu03}]\label{defn-loc-S}
A random variable $U$ and its distribution $F$ are called {\it locally
  subexponential} with span $h \in \bbN \cup \{\infty\}$ if $U \in
\calL_{\rm loc}(h)$ and
\[
\lim_{k\to\infty}
{\PP(U_1+U_2 \in k+\Delta_h) \over \PP(U \in k + \Delta_h)} 
= 2.
\]
\end{app-defn}

We denote by ${\cal S}_{\rm loc}(h)$ the class of locally
subexponential distributions with span $h$. Obviously, $\calS_{\rm
  loc}(\infty)$ is equivalent to (ordinal) subexponential class
$\calS$ (see Definition~\ref{defn-subexp}). Further,
Definition~\ref{defn-loc-S} shows that $\calS_{\rm loc}(h)\subset
\calL_{\rm loc}(h)$.

\begin{app-rem}\label{rem-local-subexp}
If $U \in \calS_{\rm loc}(h)$ for some $h \in \bbN$, then $U \in
\calS_{\rm loc}(nh)$ for all $n \in \bbN$ and $U \in \calS$ (see
Remark~2 in \cite{Asmu03}).
\end{app-rem}

\begin{app-prop}\label{prop-S^{ast}-S_{loc}}
$U \in \calS^{\ast}$ if and only if $U_{\rm e} \in \calS_{\rm loc}(1)$.
\end{app-prop}

\proof The if-part is obvious. Indeed, since $\PP(U_{\rm e} = k) =
\PP(U > k)/\EE[U]$ for $k \in \bbZ_+$, it follows that if $U_{\rm e}
\in \calS_{\rm loc}(1)$, then (\ref{def-S*}) holds, i.e., $U \in
\calS^{\ast}$.

On the other hand, suppose (\ref{def-S*}) holds for $h=1$. We then
have
\[
\lim_{k\to\infty}
\sum_{l=0}^k {\PP(U_{\rm e} = k-l)\PP(U_{\rm e} = l) \over \PP(U_{\rm e} = k)}
= 2.
\]
Further $U \in \calS \subset \calL$ (see Proposition A.2 in
\cite{Masu11}) and thus
\[
\lim_{k\to\infty}{\PP(U > k+1) \over \PP(U > k)}
=
\lim_{k\to\infty}{\PP(U_{\rm e} = k+1) \over \PP(U_{\rm e} = k)} = 1.
\]
As a result, $U_{\rm e} \in \calS_{\rm loc}(1)$. \qed

\begin{app-prop}[Proposition~3 in \cite{Asmu03}]\label{appen-prop-sum}
Suppose $U \in \calS_{\rm loc}(h)$ for some $h \in \bbN
\cup\{\infty\}$ and let $U^{(j)}$ ($j\in\bbN$) denote independent
random variables in $\bbZ_+$ such that
\[
\lim_{k\to\infty}{\PP(U^{(j)} \in k+\Delta_h) 
\over \PP(U \in k + \Delta_h)} = c_j \in \bbR_+.
\]
Then for $n \in \bbN$,
\[
\lim_{k\to\infty}{\PP(U^{(1)} + U^{(2)} + \cdots + U^{(n)} \in k+\Delta_h) 
\over \PP(U \in k + \Delta_h)}=\sum_{j=1}^nc_j.
\]
Further, if $\sum_{j=1}^nc_j > 0$, then $U^{(1)} + U^{(2)} + \cdots +
U^{(n)} \in \calS_{\rm loc}(h)$.
\end{app-prop}

\begin{app-prop}\label{prop-F(k)-bound}
Let $\{F(k);k\in\bbZ_+\}$ and $\{F_j(k);k\in\bbZ_+\}$
($j=1,2,\dots,m$) denote probability mass functions.  Suppose (i) $F
\in \calS_{\rm loc}(1)$; and (ii) for
$j=1,2,\dots,m$,
\begin{equation}
\lim_{k\to\infty}{F_j(k) \over F(k)} = c_j \in \bbR_+.
\label{eq-lim_F_j-F}
\end{equation}
Then for any $\varepsilon > 0$ there
exists some $C_{\varepsilon} \in (0,\infty)$ such that
\begin{equation}
F_1^{\ast n_1} \ast F_2^{\ast n_2} 
\ast \cdots \ast F_m^{\ast n_m}(k) 
\le C_{\varepsilon} (1+\varepsilon)^{n_1+n_2+\cdots+n_m}F(k),
\label{eq-prop-F(k)-bound}
\end{equation}
for all $k > \sup\{k\in\bbZ_+; F(k) = 0\}$ and $n_1,n_2,\dots,n_m
\in \bbN$.
\end{app-prop}

\proof See Appendix~\ref{proof-prop-F(k)-bound}. \qed

\begin{app-prop}\label{appen-prop-convo1}
For $d_i \in \bbN$ ($i=0,1,2$), let $\{\vc{P}(k);k\in\bbZ_+\}$ and
$\{\vc{Q}(k);k\in\bbZ_+\}$ denote nonnegative $d_0 \times d_1$ and
$d_1 \times d_2$ matrix sequences, respectively, such that
$\vc{P} := \sum_{k=0}^{\infty}\vc{P}(k)$ and $\vc{Q} :=
\sum_{k=0}^{\infty}\vc{Q}(k)$ are finite. Suppose that for some $U \in
\calS_{\rm loc}(1)$,
\[
\lim_{k\to\infty}{\vc{P}(k) \over \PP(U=k)} 
= \widetilde{\vc{P}} \ge \vc{O},
\quad
\lim_{k\to\infty}{\vc{Q}(k) \over \PP(U=k)} 
= \widetilde{\vc{Q}} \ge \vc{O}.
\]
We then have
\[
\lim_{k\to\infty} {\vc{P} \ast \vc{Q}(k) \over \PP(U=k)} 
= \widetilde{\vc{P}} \vc{Q} 
+ \vc{P} \widetilde{\vc{Q}}.
\]
\end{app-prop}

\proof This proposition can be proved
in the same way as Proposition A.3 in \cite{Masu11}, 
and thus the proof is omitted. \qed

\begin{app-prop}\label{appen-prop-convo2}
Let $\{\vc{W}(k);k\in\bbZ_+\}$ denote a sequence of (finite
dimensional) nonnegative square matrices such that
$\sum_{n=0}^{\infty}\vc{W}^n = (\vc{I} - \vc{W})^{-1} < \infty$, where
$\vc{W}= \sum_{k=0}^{\infty}\vc{W}(k)$. If there exists some $U \in
\calS_{\rm loc}(1)$ such that
\[
\lim_{k\to\infty}{\vc{W}(k) \over \PP(U=k)} 
= \widetilde{\vc{W}} \ge \vc{O},
\]
then
\[
\lim_{k\to\infty} {\sum_{n=0}^{\infty}\vc{W}^{\ast n}(k) 
\over \PP(U=k)} 
= (\vc{I} - \vc{W})^{-1}\widetilde{\vc{W}}
(\vc{I} - \vc{W})^{-1}.
\]
\end{app-prop}

\proof Using Proposition~\ref{appen-prop-convo1}, we can readily
prove, by induction, that
\begin{equation}
\lim_{k \to \infty} 
{ \vc{W}^{\ast n}(k) \over \PP (U=k) } 
= \sum_{l=0}^{n-1} \vc{W}^l \widetilde{\vc{W}} 
\vc{W}^{n-l-1}.
\label{lim-W^{*n}(k)}
\end{equation}
Further it follows from Proposition~\ref{prop-F(k)-bound} that for any
$\varepsilon > 0$ there exist some $k_0 \in \bbZ_+$ and some
$C_{\varepsilon} \in (0,\infty)$ such that for all $k \ge k_0$ and $n
\in \bbN$,
\[
{[\vc{W}^{\ast n}(k)]_{i,j} \over \PP(U=k)}
\le C_{\varepsilon} (1+\varepsilon)^n [\vc{W}^n]_{i,j}.
\]
Note here that $\rsp(\vc{W}) < 1$ and thus
$\sum_{n=1}^{\infty}(1+\varepsilon)^n \vc{W}^n < \infty$ for any
sufficiently small $\varepsilon > 0$. As a result, using the dominated
convergence theorem and (\ref{lim-W^{*n}(k)}), we obtain
\begin{eqnarray*}
\lim_{k \to \infty} 
{ \sum_{n = 0}^{\infty} \vc{W}^{\ast n}(k) \over \PP (U=k) } 
&=&  \lim_{k \to \infty} 
{ \vc{W}^{\ast 0}(k) \over \PP (U=k) }
+ \sum_{n = 1}^{\infty} \lim_{k \to \infty} 
{ \vc{W}^{\ast n}(k) \over \PP (U=k) }
\nonumber
\\
&=& \sum_{n = 1}^{\infty} \sum_{l=0}^{n-1} \vc{W}^l \widetilde{\vc{W}} 
\vc{W}^{n-l-1}
\nonumber
\\
&=&
(\vc{I} - \vc{W})^{-1} \widetilde{\vc{W}} (\vc{I} - \vc{W})^{-1}.
\end{eqnarray*}
\qed

\section{Proofs}

\subsection{Proof of Proposition~\ref{prop-Phi}}\label{proof-prop-Phi}

Equation (\ref{eq-RG-1}) yields
\[
\det(\vc{I} - \vc{A})
= \det(\vc{I} - \vc{R}) \det(\vc{I} - \vc{\Phi}(0)) \det(\vc{I} - \vc{G}).
\]
It thus follows from $\rsp(\vc{A}) < 1$ that
\begin{equation}
\det(\vc{I} - \vc{G}) \neq 0, \quad \det(\vc{I} - \vc{R}) \neq 0.
\label{add-eqn-11}
\end{equation}
Note here that by definition,
\[
\sum_{k=1}^N \sum_{j\in\bbM}[\vc{G}(k)]_{i,j}
= \PP(T_{<N} < \infty \mid X_0 = N, S_0 = i),\quad \mbox{for all}~ N \in \bbN,
\]
which shows that $\vc{G}\vc{e} \le \vc{e}$ and thus $\rsp(\vc{G}) \le 1$
(see Theorem~8.1.22 in \cite{Horn90}). Further, $\rsp(\vc{R}) \le 1$ due
to the duality of the $R$- and $G$-matrices (see \cite{Zhao99}).
Therefore, it follows from Theorem~8.3.1 in \cite{Horn90} and
(\ref{add-eqn-11}) that (i) $\rsp(\vc{G}) < 1$ and (ii) $\rsp(\vc{R}) <
1$.

Finally, we prove (iii). From (\ref{eqn-Phi(-l)}), we have
\[
\vc{\Phi}(-k) \ge \vc{O},
\quad 
\vc{0} \le \sum_{l=0}^{k-1}\vc{\Phi}(-l)\vc{e} \le \vc{e},
\quad \mbox{for all}~k \in \bbN,
\]
which implies that $\rsp(\sum_{l=0}^{\infty} \vc{\Phi}(-l)) \le 1$ (see
Theorem~8.1.22 in \cite{Horn90}). Thus it suffices to prove that
$\sum_{l=0}^{\infty} \vc{\Phi}(-l)$ does not have the eigenvalue
{\it one}. Indeed, (\ref{def-G(k)}) yields
\[
(\vc{I} - \vc{\Phi}(0)) (\vc{I} - \vc{G}) 
= \vc{I} - \sum_{l=0}^{\infty} \vc{\Phi}(-l).
\]
Therefore we have
$\det(\vc{I} - \sum_{l=0}^{\infty} \vc{\Phi}(-l)) \neq 0$ because
$\vc{I} - \vc{\Phi}(0)$ is nonsingular and $\rsp(\vc{G})<1$.

\subsection{Proof of Proposition~\ref{prop-structure-G}}
\label{proof-prop-structure-G}

We prove this proposition by reduction to absurdity. To do so, we
suppose either (i) $\vc{G}$ is strictly lower triangular, or (ii)
$\vc{G}$ takes a form such that
\begin{equation}
\vc{G} 
=
\left(
\begin{array}{ccc}
\vc{G}_1     & \vc{O} &       \vc{O} \\
\vc{G}_{2,1} & \vc{G}_2     & \vc{O} \\
\vc{G}_{3,1} & \vc{G}_{3,2} & \vc{G}_3
\end{array}
\right),
\label{eq-G-form-add1}
\end{equation}
where $\vc{G}_i$ $(i= 1,2)$ is irreducible and $\vc{G}_2$ can be equal
to $\vc{G}_{\rm T}$ (in that case, the last block row and column
vanish).  If (i) is true, then $\vc{G}$ is a nilpotent matrix, which
is inconsistent with $\delta(\vc{G}) = 1$.

In what follows, we consider case (ii). For simplicity, we partition
the phase set $\bbM$ into subsets $\bbM_1$, $\bbM_2$ and $\bbM_3$
corresponding to $\vc{G}_1$, $\vc{G}_2$ and $\vc{G}_3$,
respectively. Further we write $(k,i)
\stackrel{\scriptscriptstyle\neq(0,\ast)}{\longrightarrow}(l,j)$
($k,l\in\bbN; i,j \in \bbM$) when state $(l,j)$ can be reached from
state $(k,i)$ avoiding level zero.

Let $\vc{G}_2(k)$ denote a submatrix of $\vc{G}(k)$ such that
$\sum_{k=1}^{\infty}\vc{G}_2(k) = \vc{G}_2$. The irreducibility of
$\vc{G}_2$ shows that $\sum_{k=1}^{K_G}\vc{G}_2(k)$ is irreducible for
some $K_G \in \bbN$. Thus for any $i_2 \in \bbM_2$, there exists some
$(k_2',i_2') \in \bbN \times \bbM_2$ such that
\[
(k_2',i_2')\stackrel{\scriptscriptstyle\neq(0,\ast)}{\longrightarrow}
(1,i_2).
\]
Similarly, $\sum_{k=-K_A}^{\infty}\vc{A}(k)$ is irreducible for some
$K_A \in \bbN$ due to the irreducibility of $\vc{A}$, and thus there
exists some $(k_1,i_1) \in \bbN \times \bbM_1$ such that
\[
(k_1,i_1)\stackrel{\scriptscriptstyle\neq(0,\ast)}{\longrightarrow}
(k_2',i_2').
\]
As a result,
\[
(k_1,i_1) \stackrel{\scriptscriptstyle\neq(0,\ast)}{\longrightarrow}
(k_2',i_2')
\stackrel{\scriptscriptstyle\neq(0,\ast)}{\longrightarrow}
(1,i_2),
\quad i_1 \in \bbM_1,~i_2,i_2' \in \bbM_2,
\] 
which contradicts to the structure of $\vc{G}$ shown in 
 (\ref{eq-G-form-add1}).

\subsection{Proof of Proposition~\ref{prop-period-G-a}}
\label{proof-prop-period-G-a}

From Theorem~8.1.18 in \cite{Horn90}, we have 
\begin{equation}
\rsp(\widehat{\vc{R}}(\omega)) \le \delta(\widehat{\vc{R}}(1)) < 1,
\label{add-eqn-09}
\end{equation}
where the second inequality is due to the positive-recurrence of
$\vc{T}$ (see Theorem~3.4 in \cite{Zhao98}). It follows from
(\ref{eq-RG-1}), (\ref{add-eqn-09}) and $\rsp(\vc{\Phi}(0)) < 1$ that
\[
\det ( \vc{I} - \widehat{\vc{A}}(\omega) ) = 0
\iff
\det ( \vc{I} - \widehat{\vc{G}}(\omega) ) = 0.
\]
Note here that $\rsp(\widehat{\vc{A}}(\omega)) \le
\delta(\widehat{\vc{A}}(1)) = 1$ and
$\rsp(\widehat{\vc{G}}(\omega)) \le \delta(\widehat{\vc{G}}(1)) =
1$ (see Theorem~8.1.18 in \cite{Horn90}). Thus
\begin{eqnarray*}
\det (\vc{I} - \widehat{\vc{A}}(\omega)) = 0 
&\iff& \delta(\widehat{\vc{A}}(\omega)) = 1,
\\
\det (\vc{I} - \widehat{\vc{G}}(\omega)) = 0 
&\iff& \delta(\widehat{\vc{G}}(\omega)) = 1.
\end{eqnarray*}
As a result, $ \delta(\widehat{\vc{G}}(\omega)) = 1$ if and only if $
\delta(\widehat{\vc{A}}(\omega)) = 1$.  Finally, the statement (i)
follows from (\ref{def-tau_G}) and Proposition~\ref{prop-MAdP}.

Since the statement (i) is proved, we readily obtain the statements
(ii) and (iii) by applying Theorem~B.1 in \cite{Kimu10} to the MAdP
$\{(\breve{X}^{(G)}_n,\breve{S}^{(G)}_n)\}$ and using
(\ref{add-eqn-04}). Further, the statement (iv) is an immediate
consequence of (\ref{add-eqn-04}) and Lemma~B.3 in \cite{Kimu10}.

\subsection{Proof of Proposition~\ref{lem-adjG}}\label{proof-lem-adjG}

Since $\vc{A}$ is stochastic, it follows from
Propositions~\ref{prop-delta-Gamma_A-1} and \ref{prop-period-G-a} that
for $\nu=0,1,\dots,\tau-1$,
\begin{align}
\delta(\widehat{\vc{A}}(\omega_{\tau}^{\nu})) 
&= \delta(\widehat{\vc{A}}(1)) = 1, &
\delta(\widehat{\vc{G}}(\omega_{\tau}^{\nu})) 
&= \delta(\widehat{\vc{G}}(1)) = 1,
\label{add-eqn-38}
\\
\vc{\mu}(\omega_{\tau}^{\nu}) 
&= \vc{\pi}\vc{\Delta}_M(\omega_{\tau}^{\nu})^{-1}, &
\vc{v}(\omega_{\tau}^{\nu}) &= \vc{\Delta}_M(\omega_{\tau}^{\nu})\vc{e},
\label{add-eqn-37}
\end{align}
where we use $\vc{\mu}(1) = \vc{\pi}$ and $\vc{v}(1) =
\vc{e}$. Therefore, (\ref{eq-right-G}) and the second equation in
(\ref{add-eqn-37}) yield $\vc{v}(\omega_{\tau}^{\nu}) =
\vc{\Delta}_M(\omega_{\tau}^{\nu})\vc{e} =
\vc{y}(\omega_{\tau}^{\nu})$.

We now define
$\widetilde{\vc{\psi}}(\omega_{\tau}^{\nu})$ as
\begin{align}
\widetilde{\vc{\psi}}(\omega_{\tau}^{\nu})&= {\vc{\mu}(\omega_{\tau}^{\nu}) 
(\vc{I} - \widehat{\vc{R}}(\omega_{\tau}^{\nu})) 
( \vc{I} - \vc{\Phi}(0))
\over
\vc{\mu}(\omega_{\tau}^{\nu}) (\vc{I} - \widehat{\vc{R}}(\omega_{\tau}^{\nu})) 
( \vc{I} - \vc{\Phi}(0))\vc{v}(\omega_{\tau}^{\nu}) },
& \nu &= 0,1,\dots,\tau-1.
\label{eqn-psi(omega_{tau}^{nu})}
\end{align}
It can be shown that 
$\widetilde{\vc{\psi}}(\omega_{\tau}^{\nu}) =
\vc{\psi}(\omega_{\tau}^{\nu})$, whose proof is given later.
From
(\ref{eq-RG-1}), we have
\[
\vc{I} - \vc{G}(\omega_{\tau}^{\nu})
= (\vc{I} - \vc{\Phi}(0))^{-1}
(\vc{I} - \widehat{\vc{R}}(\omega_{\tau}^{\nu}))^{-1}
(\vc{I} - \widehat{\vc{A}}(\omega_{\tau}^{\nu})).
\]
Pre-multiplying (resp.~post-multiplying) the above equation by
$\widetilde{\vc{\psi}}(\omega_{\tau}^{\nu})$
(resp.~$\vc{v}(\omega_{\tau}^{\nu})$) and using (\ref{add-eqn-38}), we
can readily verify that $\widetilde{\vc{\psi}}(\omega_{\tau}^{\nu})$
($=\vc{\psi}(\omega_{\tau}^{\nu})$) and
$\vc{v}(\omega_{\tau}^{\nu})=\vc{y}(\omega_{\tau}^{\nu})$ are the
left- and right-eigenvectors of
$\widehat{\vc{G}}(\omega_{\tau}^{\nu})$ corresponding to the eigenvalue
$\delta(\widehat{\vc{G}}(\omega_{\tau}^{\nu}))=1$.  As a result, the
statement (i) holds.

As for the statement (ii), it follows from the second equation in
(\ref{add-eqn-37}) and (\ref{eqn-psi(omega_{tau}^{nu})}) that
\[
\widetilde{\vc{\psi}}(\omega_{\tau}^{\nu})
\vc{\Delta}_M(\omega_{\tau}^{\nu}) \vc{e}
= \widetilde{\vc{\psi}}(\omega_{\tau}^{\nu})
\vc{v}(\omega_{\tau}^{\nu}) = 1.
\]
Therefore the statement (ii) can be proved in the same way as the
proof of Lemma~3.2 in \cite{Kimu10}.

In what follows, we prove $\widetilde{\vc{\psi}}(\omega_{\tau}^{\nu})
= \vc{\psi}(\omega_{\tau}^{\nu})$.  For this purpose, we first show
that
\begin{equation}
\sum_{l=0}^{\infty}(\omega_{\tau}^{\nu})^l \vc{\Phi}(l)
= \vc{\Delta}_M(\omega_{\tau}^{\nu}) \sum_{l=0}^{\infty} \vc{\Phi}(l)
\vc{\Delta}_M(\omega_{\tau}^{\nu})^{-1}.
\label{add-eqn-15}
\end{equation}
The definition of $\vc{\Phi}(l)$ ($l \in \bbZ_+$) implies
\[
[\vc{\Phi}(l)]_{i,j}
= \PP(X_{T_{\downarrow l+1}} =l+1, S_{T_{\downarrow l+1}} = j 
\mid X_0 = 1, S_0 = i),
\]
where $T_{\downarrow l+1} = \inf\{n \in \bbN; X_n = l+1 <
X_m~(m=1,2,\dots,n-1)\}$. Further (\ref{add-eqn-58}) and
(\ref{add-eqn-12}) imply that for all $n\in \bbN$, the following
probability is positive only if $l \equiv p(j)-p(i)~(\Mod~\tau)$:
\begin{eqnarray*}
\PP(X_n = l+1, S_n = j, X_m \ge 1~(m=1,2,\dots,n-1)
\mid X_0 = 1, S_0 = i).
\end{eqnarray*}
Thus $[\vc{\Phi}(l)]_{i,j} > 0$ only if $l \equiv
p(j)-p(i)~(\Mod~\tau)$, which leads to
\begin{eqnarray}
\sum_{l=0}^{\infty}z^l \vc{\Phi}(l)
&=& \vc{\Delta}_M(z) 
\vc{\Lambda}_{\Phi}(z^{\tau})
\vc{\Delta}_M(z)^{-1},
\label{eqn-widehat{Phi}(z)}
\end{eqnarray}
where $\vc{\Lambda}_{\Phi}(z)$ denotes an $M \times M$ matrix whose
$(i,j)$th element is given by
\[
[\vc{\Lambda}_{\Phi}(z)]_{i,j} 
= \sum_{\scriptstyle n \in \bbZ_+ \atop \scriptstyle n\tau+p(j)-p(i) \ge 0}
z^n [\vc{\Phi}(n\tau+p(j)-p(i))]_{i,j}.
\]
As a result, (\ref{eqn-widehat{Phi}(z)}) yields (\ref{add-eqn-15})
because $\vc{\Lambda}_{\Phi}(1) = \sum_{l=0}^{\infty}\vc{\Phi}(l)$.

We now return to the proof of
$\widetilde{\vc{\psi}}(\omega_{\tau}^{\nu}) =
\vc{\psi}(\omega_{\tau}^{\nu})$. From (\ref{def-R(k)}) and
(\ref{add-eqn-15}), we have for $\nu = 0,1,\dots,\tau-1$,
\begin{eqnarray}
\lefteqn{
(\vc{I} - \widehat{\vc{R}}(\omega_{\tau}^{\nu})) 
( \vc{I} - \vc{\Phi}(0))
}
\quad &&
\nonumber
\\
&=& \vc{I} - \sum_{l=0}^{\infty}(\omega_{\tau}^{\nu})^l \vc{\Phi}(l)
\nonumber
\\
&=& \vc{\Delta}_M(\omega_{\tau}^{\nu}) 
\left(\vc{I} - \sum_{l=0}^{\infty} \vc{\Phi}(l) \right)
\vc{\Delta}_M(\omega_{\tau}^{\nu})^{-1}
\nonumber
\\
&=& \vc{\Delta}_M(\omega_{\tau}^{\nu}) 
(\vc{I} - \vc{R}) ( \vc{I} - \vc{\Phi}(0)) 
\vc{\Delta}_M(\omega_{\tau}^{\nu})^{-1},
\label{add-eqn-17}
\end{eqnarray}
where the last equality follows from the first equality with $\nu=0$.
Substituting (\ref{add-eqn-37}) and (\ref{add-eqn-17}) into
(\ref{eqn-psi(omega_{tau}^{nu})}) yields
\[
\widetilde{\vc{\psi}}(\omega_{\tau}^{\nu}) 
= {\vc{\pi}(\vc{I} - \vc{R})(\vc{I} - \vc{\Phi}(0))
\over
\vc{\pi}(\vc{I} - \vc{R})(\vc{I} - \vc{\Phi}(0))\vc{e} }
\vc{\Delta}_M(\omega_{\tau}^{\nu})^{-1} 
= \vc{\psi}(\omega_{\tau}^{\nu}).
\]

\subsection{Proof of Lemma~\ref{lem-lim-D}}\label{proof-lem-lim-D}

From (\ref{defn-hat{L}(z)}), we have
\begin{equation}
\widehat{\vc{L}}(1/z)
= {\adj(\vc{I} - \widehat{\vc{G}}(1/z)) 
\over \det(\vc{I} - \widehat{\vc{G}}(1/z))} - \vc{I}.
\label{eq-D^hat-G-02}
\end{equation}
Note here that
\begin{align*}
\left| [\widehat{\vc{G}}(1/z)]_{i,j} \right|
&= \left| \sum_{k=1}^{\infty}z^k[\vc{G}(k)]_{i,j} \right|
\le [\vc{G}]_{i,j}, & i,j &\in \bbM,~|z| \le 1,
\\
\rsp(\widehat{\vc{G}}(1/z)) &< \rsp(\vc{G}) = 1, & |z| & < 1.
\end{align*}
It then follows from Proposition~\ref{prop-period-G-a} that
$\{\omega_{\tau}^{\nu};\nu=0,1,\dots,\tau-1\}$ are the simple
minimum-modulus poles of $\widehat{\vc{L}}(1/z)$.
Therefore applying Theorem~A.1 in \cite{Kimu10} to
(\ref{eq-D^hat-G-02}), we obtain
\begin{equation}
\vc{L}(k) 
=
\sum_{ \nu = 0}^{\tau-1} 
{1 \over (\omega_{\tau}^{\nu})^k}
\lim_{z \to \omega_{\tau}^{\nu} } 
\left( 1 - { z \over \omega_{\tau}^{\nu}}  \right) 
{\adj (\vc{I} - \widehat{\vc{G}}(1/z) ) 
\over \det (\vc{I} - \widehat{\vc{G}}(1/z) )}
+ O ( (1 + \varepsilon_0)^{-k} ) \vc{e}^{\rmt} \vc{e}, 
\label{eqn-lim-D-1} 
\end{equation}
for some $\varepsilon_0 > 0$, where $f(x) = O(g(x))$ represents
$\limsup_{x\to\infty}|f(x)/g(x)| < \infty$. Further it follows from
l'H\^{o}pital's rule and Proposition~\ref{lem-adjG} that for $\nu =
0,1,\dots,\tau-1$,
\begin{eqnarray}
\lefteqn{
\lim_{z \to \omega_{\tau}^{\nu} } \left( 1 - { z \over \omega_{\tau}^{\nu}}  \right) 
{\adj (\vc{I} - \widehat{\vc{G}}(1/z) ) 
\over \det (\vc{I} - \widehat{\vc{G}}(1/z) )}
}
\quad \nonumber
\\
&=&
\lim_{z \to \omega_{\tau}^{\nu} } 
{1-\displaystyle{ z \over \omega_{\tau}^{\nu}} 
\over 1 - \delta( \widehat{\vc{G}}(1/z) )}
{\adj (\vc{I} -\widehat{\vc{G}}(\omega_{\tau}^{-\nu})) 
\over \displaystyle{\prod_{i=2}^{M} ( 1- \lambda_i^{(G)}(\omega_{\tau}^{-\nu}))}} 
\nonumber
\\
&=& 
{
1 \over \omega_{\tau}^{\nu} \cdot 
(\rd/\rd z) \delta( \widehat{\vc{G}}(1/z) )|_{z =\omega_{\tau}^{\nu}} }
\cdot \vc{y}(\omega_{\tau}^{-\nu} )\vc{\psi} (\omega_{\tau}^{-\nu})
\nonumber
\\
&=&
{
1 \over \omega_{\tau}^{\nu} \cdot 
(\rd/\rd z) \delta( \widehat{\vc{G}}(1/z) )|_{z =\omega_{\tau}^{\nu}} }
\cdot
\vc{\Delta}_M(\omega_{\tau}^{-\nu}) 
{ \vc{e}\vc{\psi} \over \vc{\psi} \vc{e} }
 \vc{\Delta}_M(\omega_{\tau}^{-\nu})^{-1},
\label{eqn-lim-D-add-1}
\end{eqnarray}
where the last equality is due to (\ref{eq-left-G}), (\ref{eq-right-G}) and
(\ref{defn-psi}).  Letting $y=1/z$, we have
\begin{eqnarray}
\left.
\omega_{\tau}^{\nu}
 {\rd \over \rd z} \delta( \widehat{\vc{G}}(1/z) ) 
\right|_{z = \omega_{\tau}^{\nu}}
= -
\left.{1 \over \omega_{\tau}^{\nu}}
 {\rd \over \rd y} \delta( \widehat{\vc{G}}(y)) 
\right|_{y = 1/\omega_{\tau}^{\nu}}
= -
\left. {\rd \over \rd y} \delta( \widehat{\vc{G}}(y)) \right|_{y = 1},
\nonumber
\\
\label{add-eqn-39}
\end{eqnarray}
where the second equality is due to
Proposition~\ref{prop-period-G-a}~(iv). Applying (\ref{add-eqn-39}) to
(\ref{eqn-lim-D-add-1}) yields
\begin{eqnarray}
\lefteqn{
\lim_{z \to \omega_{\tau}^{\nu} } 
\left( 1 - { z \over \omega_{\tau}^{\nu}} \right) 
{\adj( \vc{I} - \widehat{\vc{G}}(1/z) ) 
\over \det( \vc{I} - \widehat{\vc{G}}(1/z) )}
}
\quad &&
\nonumber
\\
&=& 
 {-1 \over (\rd/\rd y) \delta( \widehat{\vc{G}}(y) )|_{y=1} }
\cdot
\vc{\Delta}_M(\omega_{\tau}^{-\nu}) 
{ \vc{e}\vc{\psi} \over \vc{\psi} \vc{e} }
 \vc{\Delta}_M(\omega_{\tau}^{-\nu})^{-1}.
\label{eq-lim-adjdetG}
\end{eqnarray}

In what follows, we calculate $(\rd/\rd y) \delta( \widehat{\vc{G}}(y)
)|_{y=1}$. Taking the derivative of both sides of
(\ref{def-eigenvalues-G}) with $z=y$, letting $y =1$ and using
$\delta(\widehat{\vc{G}}(1))=1$, we have
\begin{equation}
\left.{\rd\over\rd y} \delta( \widehat{\vc{G}}(y) )\right|_{y=1}
= -
{1 \over \displaystyle{\prod_{i=2}^{M} (1 - \lambda_i^{(G)}(1))} }
\cdot
{\rd \over \rd y}\left. \det( \vc{I} - \widehat{\vc{G}}(y)) \right|_{y=1}. 
\label{eq-add-delta-G-1}
\end{equation}
Similarly, from $\det(\vc{I} - \widehat{\vc{G}}(y))= \vc{\pi}\cdot
\adj(\vc{I} - \widehat{\vc{G}}(y))(\vc{I}-\widehat{\vc{G}}(y)) \cdot
\vc{e}$, we obtain
\begin{eqnarray}
{\rd \over \rd y}\left. 
\det( \vc{I} - \widehat{\vc{G}}(y)) \right|_{y=1}
&=& \vc{\pi} \cdot \adj(\vc{I} - \vc{G}) 
\sum_{k=1}^{\infty}k\vc{G}(k)\vc{e},
\label{add-eqn-05a}
\end{eqnarray}
where we use $\vc{G}\vc{e}=\vc{e}$.  Note here that
Proposition~\ref{lem-adjG} and (\ref{defn-psi}) imply
\[
\adj(\vc{I} - \vc{G})
= { \vc{e}\vc{\psi} \over \vc{\psi} \vc{e} }
\cdot \prod_{i=2}^{M} (1 - \lambda_i^{(G)}(1)).
\]
It thus follows from (\ref{add-eqn-05a}) and Lemma~\ref{lem-sigma}
that
\begin{eqnarray}
{\rd \over \rd y}\left. 
\det( \vc{I} - \widehat{\vc{G}}(y)) \right|_{y=1}
&=& { \vc{\psi} \over \vc{\psi} \vc{e} }
\sum_{k=1}^{\infty}k\vc{G}(k)\vc{e}
\cdot \prod_{i=2}^{M} (1 - \lambda_i^{(G)}(1))
\nonumber
\\
&=& { 1 \over \vc{\psi} \vc{e} }
\cdot \prod_{i=2}^{M} (1 - \lambda_i^{(G)}(1)),
\label{add-eqn-05b}
\end{eqnarray}
where the second equality is due to
$\vc{\psi}\sum_{k=1}^{\infty}k\vc{G}(k)\vc{e} = 1$ (see
(\ref{eqn-sigma}) and (\ref{defn-psi})).  Further substituting
(\ref{add-eqn-05b}) into (\ref{eq-add-delta-G-1}) yields
\[
\left.{\rd\over\rd y} \delta( \widehat{\vc{G}}(y) )\right|_{y=1}
=
-{1 \over \vc{\psi} \vc{e} },
\]
from which and (\ref{eq-lim-adjdetG}), we have
\begin{eqnarray}
\lim_{z \to \omega_{\tau}^{\nu} } 
\left( 1 - { z \over \omega_{\tau}^{\nu}}  \right) 
{\adj (\vc{I} - \widehat{\vc{G}}(1/z) )
\over \det( \vc{I} - \widehat{\vc{G}}(1/z) )}
&=&
\vc{\Delta}_M(\omega_{\tau}^{-\nu}) \vc{e}\vc{\psi}
\vc{\Delta}_M(\omega_{\tau}^{-\nu})^{-1}. \qquad
 \label{eq-lim-adjdetG-3}
\end{eqnarray}
Finally, we have (\ref{eq-lim-D}) by substituting
(\ref{eq-lim-adjdetG-3}) into (\ref{eqn-lim-D-1}) and letting $k = n
\tau + l$. 

\subsection{Proof of Lemma~\ref{lem-AD-BD}}\label{proof-lem-AD=BD}

Equations (\ref{eqn-lim-sum-D(ntau+l)}) and (\ref{eqn-lim-A(k)e}) show
that for any $\varepsilon >0$ there exists some $m_{\ast}:= m_{\ast}
(\varepsilon) \in \bbN$ such that for all $m \ge m_{\ast}$ and $l =
0,1,\dots,\tau-1$,
\begin{eqnarray}
&&
\vc{e}(\tau\vc{\psi} - \varepsilon\vc{e}^{\rmt})
\le
\sum_{l=0}^{\tau-1}
\vc{L} ( \lfloor m / \tau \rfloor \tau + l) 
\le \vc{e}(\tau\vc{\psi} + \varepsilon\vc{e}^{\rmt}),\qquad
\label{eq-ep-D}
\\
&&
{1 \over \EE [Y] } (\vc{c}{}_{\!A} - \varepsilon \vc{e}) 
\le
{ \overline{\vc{A}}( \lfloor m / \tau \rfloor \tau + l)\vc{e}  \over \PP(Y> m) }
\leq
{1 \over \EE [Y] } (\vc{c}{}_{\!A} + \varepsilon \vc{e}).
\label{eq-ep-barB}
\end{eqnarray}
Further since $Y_{\rm e} \in \calL$ and $\vc{L}(m) \leq
\vc{e}\vc{e}^{\rmt}$ for all $m = 1,2,\dots$, we have
\begin{eqnarray}
\lefteqn{
\limsup_{k \to \infty} \sum_{m=1}^{m_{\ast}-1} 
{ \overline{\vc{A}}(k+m) \vc{L}(m) \over \PP(Y_{\rm e} > k) }
}
\quad &&
\nonumber
\\
&\leq&
\sum_{m=1}^{m_{\ast}-1} \limsup_{k \to \infty} 
{ \overline{\vc{A}}(k+m) \vc{e}\vc{e}^{\rmt} \over \PP(Y > k+m) }
\limsup_{k \to \infty} {\PP(Y > k+m) \over  \PP(Y_{\rm e} > k+m) } 
\nonumber
\\
&& {} \quad \times
\limsup_{k \to \infty} {\PP(Y_{\rm e} > k+m) \over \PP(Y_{\rm e} > k) }
\nonumber
\\
&=&
\vc{O},
\label{eqn-lim-AD-1}
\end{eqnarray}
where the last equality follows from (\ref{eqn-lim-A(k)e}) and the
fact that $Y_{\rm e} \in \calL$ has a heavier tail than $Y$ (see
Corollary 3.3 in \cite{Sigm99}).

On the other hand, 
\begin{eqnarray}
\sum_{m=m_{\ast}}^{\infty}
{ \overline{\vc{A}} (k+m) \vc{L}(m) 
\over \PP(Y_{\rm e} > k) }
&\le&
\sum_{m'=  \lfloor m_{\ast}/\tau \rfloor}^{\infty}
\sum_{l=0}^{\tau-1}
{ \overline{\vc{A}}(k + m'\tau + l) \vc{L}(m'\tau + l) 
\over \PP(Y_{\rm e} > k) }
\nonumber
\\
&\le& 
\sum_{m'=  \lfloor m_{\ast}/\tau \rfloor}^{\infty}
{ \overline{\vc{A}}(k + m'\tau) 
\over \PP(Y_{\rm e} > k) }
\sum_{l=0}^{\tau-1}
\vc{L}(m'\tau + l)
\nonumber
\\
&\le&
\sum_{m'=  \lfloor m_{\ast}/\tau \rfloor}^{\infty}
{ \overline{\vc{A}}(k + m'\tau) \vc{e}
\over \PP(Y_{\rm e} > k) }
(\tau\vc{\psi} + \varepsilon \vc{e}^{\rmt}),
\label{add-eqn-07}
\end{eqnarray}
where the second inequality holds because $\{\overline{\vc{A}}
(k);k\in\bbZ_+\}$ is nonincreasing, and where the last inequality is
due to (\ref{eq-ep-D}).  Note here that (\ref{eqn-lim-A(k)e}) implies
for all sufficiently large $k$,
\begin{eqnarray*}
\lefteqn{
\sum_{m'=  \lfloor m_{\ast}/\tau \rfloor}^{\infty}
{ \overline{\vc{A}}(k + m'\tau) \vc{e}
\over \PP(Y_{\rm e} > k) }
}
\qquad &&
\nonumber
\\
&\le& 
(\vc{c}{}_{\!A} + \varepsilon \vc{e})
\cdot
{1 \over \EE [Y] } \sum_{m'=  \lfloor m_{\ast}/\tau \rfloor}^{\infty}
{ \PP(Y > k + m'\tau)
\over \PP(Y_{\rm e} > k) },
\end{eqnarray*}
from which and Proposition \ref{add-appen-prop} it follows that
\begin{equation}
\limsup_{k\to\infty}
\sum_{m'=  \lfloor m_{\ast}/\tau \rfloor}^{\infty}
{ \overline{\vc{A}}(k + m'\tau) \vc{e}
\over \PP(Y_{\rm e} > k) }
\le
{\vc{c}{}_{\!A} + \varepsilon \vc{e} \over \tau}.
\label{add-eqn-40}
\end{equation}
Combining (\ref{add-eqn-07}) and (\ref{add-eqn-40}) and letting
$\varepsilon \downarrow 0$ yield
\begin{eqnarray}
\limsup_{k\to\infty}
\sum_{m=m_{\ast}}^{\infty}
{ \overline{\vc{A}} (k+m) \vc{L}(m) 
\over \PP(Y_{\rm e} > k) }
&\le& \vc{c}{}_{\!A}\vc{\psi}.
\label{add-eqn-23}
\end{eqnarray}
As a result, from (\ref{eqn-lim-AD-1}) and (\ref{add-eqn-23}), we have
\begin{eqnarray}
\limsup_{k \to \infty}
\sum_{m=1}^{\infty} 
{\overline{\vc{A}}(k+m) \vc{L}(m) \over \PP(Y_{\rm e} > k) }
&\leq& 
\vc{c}{}_{\!A}\vc{\psi}.
\label{add-eqn-13}
\end{eqnarray}

Next we consider the lower limit. It follows from (\ref{eq-ep-D}) and
(\ref{eq-ep-barB}) that
\begin{eqnarray}
\sum_{m=1}^{\infty}
{ \overline{\vc{A}} (k+m) \vc{L}(m) 
\over \PP(Y_{\rm e} > k) }
&\ge&
\sum_{m=m_{\ast}}^{\infty}
{ \overline{\vc{A}} (k+m) \vc{L}(m) 
\over \PP(Y_{\rm e} > k) }
\nonumber
\\
&\ge&
\sum_{m'=  \lfloor m_{\ast}/\tau \rfloor + 1}^{\infty}
\sum_{l=0}^{\tau-1}
{ \overline{\vc{A}}(k + m'\tau + l) \vc{L}(m'\tau + l) 
\over \PP(Y_{\rm e} > k) }
\nonumber
\\
&\ge& 
\sum_{m'=  \lfloor m_{\ast}/\tau \rfloor + 1}^{\infty}
{ \overline{\vc{A}}(k + m'\tau + \tau) 
\over \PP(Y_{\rm e} > k) }
\sum_{l=0}^{\tau-1}
\vc{L}(m'\tau + l)
\nonumber
\\
&\ge& 
\sum_{m'=  \lfloor m_{\ast}/\tau \rfloor+2}^{\infty}
{ \overline{\vc{A}}(k + m'\tau) \vc{e}
\over \PP(Y_{\rm e} > k) }
(\tau\vc{\psi} - \varepsilon \vc{e}^{\rmt}),
\label{add-eqn-41}
\end{eqnarray}
where the third inequality requires the fact that
$\{\overline{\vc{A}}(k)\}$ is nonincreasing. Further the following can
be shown in a very similar way to (\ref{add-eqn-40}):
\[
\liminf_{k\to\infty}
\sum_{m'=  \lfloor m_{\ast}/\tau \rfloor + 2}^{\infty}
{ \overline{\vc{A}}(k + m'\tau) \vc{e}
\over \PP(Y_{\rm e} > k) }
\ge
{\vc{c}{}_{\!A} - \varepsilon \vc{e} \over \tau}.
\]
Combining this with (\ref{add-eqn-41}) and letting $\varepsilon
\downarrow 0$ yield
\begin{equation}
\liminf_{k \to \infty}
\sum_{m=1}^{\infty} {\overline{\vc{A}}(k+m) \vc{L}(m) \over \PP(Y_{\rm e} > k) }
\geq
\vc{c}{}_{\!A}\vc{\psi}.
\label{eq-R-A-add4}
\end{equation}
Finally, (\ref{limit-AD}) follows from (\ref{add-eqn-13}),
(\ref{eq-R-A-add4}) and (\ref{defn-psi}). Equation (\ref{limit-BD})
can be proved in the same way, and thus the proof is omitted.

\subsection{Proof of Lemma~\ref{sec4-lem-AD-BD}}\label{proof-sec4-lem-AD-BD}

We give the proof of (\ref{sec4-limit-AD}) only. Equation
(\ref{sec4-limit-BD}) can be proved in the same way. It follows from
(\ref{eqn-lim-D(ntau+l)}), $\vc{E}\vc{e}=\vc{e}$ and
(\ref{assu-sec4-A(k)E-B(k)E}) that for $\varepsilon > 0$ there exists
some $m_{\ast}:=m_{\ast}(\varepsilon) \in \bbN$ such that for all $m
\ge m_{\ast}$ and $l=0,1,\dots,\tau-1$,
\begin{eqnarray}
\vc{E}(\tau\vc{H}_l - \varepsilon \vc{e} \vc{e}^{t}) 
&\le&
\vc{L}(m)
\leq
\vc{E}(\tau\vc{H}_l + \varepsilon \vc{e} \vc{e}^{t}),~~
m \equiv l~(\mbox{mod}~\tau),\quad
\label{add-eqn-21a}
\\
{1 \over \EE [Y] } (\vc{C}{}_{\!A}^{\rmE} - \varepsilon \vc{e} \vc{e}^{t}) 
&\le&
{ \vc{A}(m)\vc{E}  \over \PP(Y = m) }
\leq
{1 \over \EE [Y] } (\vc{C}{}_{\!A}^{\rmE} + \varepsilon \vc{e} \vc{e}^{\rmt}).
\label{add-eqn-21b}
\end{eqnarray}
Thus from (\ref{sec4-eq-assu-AE-BE}), $\vc{L}(m) \le
\vc{E}\vc{e}\vc{e}^{\rmt}$ and $Y \in \calL$ (see
Remark~\ref{rem-L_loc}), we have
\begin{eqnarray*}
\lefteqn{
\lim_{k\to \infty} \sum_{m=1}^{m_{\ast}-1}
{ \vc{A} (k+m) \vc{L}(m) 
\over \PP(Y_{\rm e} = k) }
}
\quad &&
\nonumber
\\
&\le& \EE[Y]
\sum_{m=1}^{m_{\ast}-1}
\lim_{k\to \infty} { \vc{A} (k+m)\vc{E}\vc{e}\vc{e}^{\rmt} \over \PP(Y = k+m) }
{\PP(Y = k+m) \over \PP(Y > k)}
= \vc{O}.
\end{eqnarray*}
Using this and (\ref{add-eqn-21a}), we obtain
\begin{eqnarray}
\lefteqn{
\limsup_{k\to \infty} \sum_{m=1}^{\infty} 
{ \vc{A} (k+m) \vc{L}(m) 
\over \PP(Y_{\rm e} = k) }
}
\quad &&
\nonumber
\\
&=& \limsup_{k\to \infty} \sum_{m=m_{\ast}}^{\infty} 
{ \vc{A} (k+m) \vc{L}(m) 
\over \PP(Y_{\rm e} = k) }
\nonumber
\\
&=&
\limsup_{k\to \infty} \sum_{l = 0}^{\tau - 1} 
\underset{m \equiv l\,(\Mod\,\tau) }{\sum_{m \geq m_{\ast}}}
{ \vc{A} (k+m) \vc{L}(m) \over \PP(Y_{\rm e} = k) }
\nonumber
\\
&\le& 
\sum_{l = 0}^{\tau - 1} 
\left[\limsup_{k \to \infty} 
\underset{m \equiv l\,(\Mod\,\tau) }{\sum_{m \geq m_{\ast}}}
{ \vc{A}(k+m)\vc{E} \over \PP(Y_{\rm e} = k) }
\right] (\tau\vc{H}_l + \varepsilon \vc{e}\vc{e}^{\rmt}).
\label{sec4-eq-A-D-1}
\end{eqnarray}
Further it follows from (\ref{add-eqn-21b}) and
Proposition~\ref{add-appen-prop2} that
\begin{eqnarray}
\lefteqn{
\limsup_{k \to \infty} 
\underset{m \equiv l\,(\Mod\,\tau) }{\sum_{m \geq m_{\ast}}}
{\vc{A}(k+m)\vc{E} \over \PP(Y_{\rme} = k) }
}
\quad
\nonumber
\\
&\le&
{\vc{C}{}_{\!A}^{\rmE} + \varepsilon \vc{e} \vc{e}^{t} \over \EE[Y]}
\limsup_{k \to \infty}
\underset{m \equiv l\,(\Mod\,\tau) }{\sum_{m \geq m_{\ast}}}
{\PP(Y = k+m) \over \PP(Y_{\rme} = k) }
\nonumber
\\
&=& {\vc{C}{}_{\!A}^{\rmE} + \varepsilon \vc{e} \vc{e}^{t} \over \tau }.
\label{sec4-eq-A-calA}
\end{eqnarray}
Substituting (\ref{sec4-eq-A-calA}) into (\ref{sec4-eq-A-D-1}) and
letting $\varepsilon \downarrow 0$, we obtain
\[
\limsup_{k\to \infty} \sum_{m=1}^{\infty} 
{ \vc{A} (k+m) \vc{L}(m) 
\over \PP(Y_{\rm e} = k) }
\le \vc{C}{}_{\!A}^{\rmE}\sum_{l = 0}^{\tau - 1}\vc{H}_l
= \vc{C}{}_{\!A}^{\rmE}\vc{e}\vc{\psi},
\]
where we use (\ref{sum-H_l}) in the last equality. 
Similarly, we can show that
\[
\liminf_{k\to \infty} \sum_{m=1}^{\infty} 
{ \vc{A} (k+m) \vc{L}(m) \over \PP(Y_{\rm e} = k) }
\ge \vc{C}{}_{\!A}^{\rmE}\vc{e}\vc{\psi}.
\]
As a result, 
\[
\lim_{k\to \infty} \sum_{m=1}^{\infty} 
{ \vc{A} (k+m) \vc{L}(m) \over \PP(Y_{\rm e} = k) }
= \vc{C}{}_{\!A}^{\rmE}\vc{e}\vc{\psi},
\]
from which and (\ref{defn-psi}) we have (\ref{sec4-limit-AD}).

\subsection{Proof of Proposition~\ref{add-appen-prop2}}\label{proof-add-appen-prop2}

We assume that condition (i) holds. It follows from $U \in \calL_{\rm
  loc}(1)$ that for any $\varepsilon > 0$ there exists $k_0 \in \bbN$
such that for all $k \ge k_0$ and $l \in \bbZ_+$,
\[
1 - \varepsilon 
\le
{\PP(U = k + lh + \nu) \over \PP(U = k + lh)}
\le 1 + \varepsilon,\qquad  \nu=0,1,\dots,h-1.
\]
Thus for all $k \ge k_0$, we have
\[
1 - \varepsilon
\le 
{\sum_{l=l_0}^{\infty}\PP(U = k + lh + \nu) 
\over 
\sum_{l=l_0}^{\infty}\PP(U = k + lh)}
\le
1 + \varepsilon,
\qquad  \nu=0,1,\dots,h-1,
\]
which leads to
\begin{equation}
\lim_{k\to\infty}
{\sum_{l=l_0}^{\infty}\PP(U = k + lh + \nu) 
\over 
\sum_{l=l_0}^{\infty}\PP(U = k + lh)} = 1,
\qquad  \nu=0,1,\dots,h-1.
\label{add-eqn-23'}
\end{equation}
Therefore (\ref{add-eqn-23'}) yields for $\nu=0,1,\dots,h-1$,
\begin{eqnarray}
\lefteqn{
\lim_{k\to\infty}
{\sum_{l=l_0}^{\infty}\PP(U = k + lh + \nu) \over \PP(U > k + l_0h-1)}
}
\quad &&
\nonumber
\\
&=&
\lim_{k\to\infty}
{\sum_{l=l_0}^{\infty}\PP(U = k + lh + \nu) 
\over \sum_{m=l_0h}^{\infty}\PP(U = k + m)}
\nonumber
\\
&=&
\lim_{k\to\infty}
{\sum_{l=l_0}^{\infty}\PP(U = k + lh) 
\over \sum_{j=0}^{h-1}\sum_{l=l_0}^{\infty}\PP(U = k + lh + j)}
\nonumber
\\
&& {} \qquad \times 
{\sum_{l=l_0}^{\infty}\PP(U = k + lh + \nu) 
\over \sum_{l=l_0}^{\infty}\PP(U = k + lh)}
 = {1 \over h}.
\label{add-eqn-24'}
\end{eqnarray}
Note here that if $U \in \calL_{\rm loc}(1)$, then $U \in \calL$ and
thus $\lim_{k\to\infty}\PP(U > k + l_0h-1) / \PP(U > k)=1$. As a
result, (\ref{add-eqn-24'}) implies (\ref{add-eqn-24}).

Next we assume that condition (ii) holds.  It then follows that for
all sufficiently large $k$,
\begin{equation}
\sum_{l=l_0}^{\infty}\PP(U = k + lh) 
\ge \sum_{l=l_0}^{\infty}\PP(U = k + lh+ j),
\qquad j \in \bbZ_+.
\label{add-eqn-42}
\end{equation}
Thus
for any
fixed (possibly negative) integer $i$,
\begin{eqnarray*}
\lefteqn{
\lim_{k\to\infty}
{\PP(U = k + l_0h + i)
\over h\sum_{l=l_0}^{\infty}\PP(U = k + lh)}
}
\quad &&
\nonumber
\\
&\le& \lim_{k\to\infty}
{\PP(U = k + l_0h + i)
\over \sum_{j=0}^{h-1}\sum_{l=l_0}^{\infty}\PP(U = k + lh+ j)}
\nonumber
\\
&=& \lim_{k\to\infty}
{\PP(U > k + l_0h + i-1) - \PP(U > k + l_0h + i) \over \PP(U > k+l_0h-1)}
= 0,
\end{eqnarray*}
which implies that
\begin{equation}
\lim_{k\to\infty}
{\PP(U = k + l_0h + i)
\over \sum_{l=l_0}^{\infty}\PP(U = k + lh)} = 0.
\label{add-eqn-52}
\end{equation}
Further (\ref{add-eqn-42}) yields for all sufficiently large $k$, 
\begin{eqnarray*}
1
&\ge&
{\sum_{l=l_0}^{\infty}\PP(U = k + lh + \nu) 
\over 
\sum_{l=l_0}^{\infty}\PP(U = k + lh)}
\nonumber
\\
&\ge&
1 -
{\PP(U = k + l_0h) 
\over 
\sum_{l=l_0}^{\infty}\PP(U = k + lh)}, \qquad \nu=0,1,\dots,h-1,
\end{eqnarray*}
from which and (\ref{add-eqn-52}) it follows that (\ref{add-eqn-23'})
holds for $\nu=0,1,\dots,h-1$. Therefore we can prove
(\ref{add-eqn-24}) in the same way as the case of condition (i).

\subsection{Proof of Proposition~\ref{prop-F(k)-bound}}\label{proof-prop-F(k)-bound}

The techniques for the proof are based on Lemma 4.2 in \cite{Asmu94}
and Lemma~10 in \cite{Jele98}, though some modifications are
required. For the reader's convenience, we provide a complete proof of
this proposition.

We first prove the statement under an additional condition that $c_j >
0$ for all $j=1,2,\dots,m$, and then remove the condition.

Let $C = \max \{1, c_1, \dots, c_{m} \}$, $d_0 = 1$ and $d_j = c_j /C
\le 1$ for $j=1,2,\dots, m$. Let $F_0(k)$ ($k\in\bbZ_+$) denote a
probability mass function such that $F_0(k)=CF(k)$ for all
sufficiently large $k \ge k_0$, where $k_0$ is a positive integer such
that $F(k) > 0$ for all $k \ge k_0$ (see Definitions~\ref{defn-loc-L}
and \ref{defn-loc-S}).

From (\ref{eq-lim_F_j-F}), we have
\begin{equation}
\lim_{k\to\infty}{F_j(k) \over F_0(k)} = d_j \le 1, \qquad j=0,1,\dots,m.
\label{eq-lim_F_j-F-02}
\end{equation}
Further since $F_j \in \calS_{\rm loc}(1) \subset \calL_{\rm loc}(1)$
(see Proposition~\ref{appen-prop-sum}),
\begin{eqnarray}
\lim_{n\to\infty}
\lim_{k\to\infty}
{\sum_{l=0}^n  F_i(l) F_j(k - l) \over F_j (k) }
&=& \lim_{n\to\infty}\sum_{l=0}^n F_i(l) = 1,
\label{add-eqn-51}
\\
\lim_{k\to\infty}
{F_i \ast F_j(k) \over F_0(k) }
&=& d_i+d_j, 
\label{add-eqn-53}
\end{eqnarray}
for all $i,j =0,1,\dots,m$.  Thus any $\varepsilon >0$, there exist
some positive integers $k'$ and $k''$ such that $k'' > 2k' \ge 2k_0$,
$F_0(k) = C F(k) \leq 1$ for all $k \ge k'$ and for all
$i,j=0,1,\dots,m$,
\begin{align}
{F_0(k+1) \over F_0(k)} & \ge 1 - \varepsilon, & \forall k &\ge k',
\label{eq-add-asum0}
\\
d_j - {\varepsilon \over 8} \le
{F_j(k) \over F_0(k)}
& \le 1 + {\varepsilon \over 2}, 
& \forall k &\ge k',
\label{eq-add-asum1}
\\
{\sum_{l=0}^{k'-1} F_i(l) F_j(k - l)
\over F_j (k) }
& \ge 1 - {\varepsilon \over 8 d_j},
& \forall k &\ge k'',
\label{eq-add-asum2}
\\
F_i \ast F_j(k) 
& \le (d_i + d_j + \varepsilon/4)F_0(k), 
& \forall k &\ge k''.
\label{eq-add-asum3}
\end{align}
Note here that (\ref{eq-add-asum0}), (\ref{eq-add-asum1}),
(\ref{eq-add-asum2}) and (\ref{eq-add-asum3}) follow from $F_0 \in
\calL_{\rm loc}(1)$, (\ref{eq-lim_F_j-F-02}), (\ref{add-eqn-51}) and
(\ref{add-eqn-53}), respectively.

We now show (\ref{eq-prop-F(k)-bound}) for the convolution of two mass
functions $F_i$ and $F_j$ ($i,j=0,1,\dots,m$). Note that
\begin{equation}
F_i \ast F_j(k)
= \sum_{l=0}^{k-k'} F_i (k-l) F_j(l) + \sum_{l=0}^{k'-1} F_i(l) F_j(k-l).
\label{add-eqn-54}
\end{equation}
It then follows from (\ref{eq-add-asum1}),
(\ref{eq-add-asum2}) and (\ref{eq-add-asum3}) that for $k \ge k'' > 2
k'$,
\begin{eqnarray}
\sum_{l=0}^{k-k'} F_i (k-l) F_j(l) 
&=& F_i \ast F_j(k) - \sum_{l=0}^{k'-1} F_i(l) F_j(k-l)
\nonumber
\\
&\le& \left( d_i + d_j +{\varepsilon \over 4} \right)F_0(k) 
- \left( 1 - {\varepsilon \over 8d_j} \right)F_j(k)
\nonumber
\\
&\le& 
\left[
\left( d_i + d_j +{\varepsilon \over 4} \right) 
- \left( 1 - {\varepsilon \over 8d_j} \right)
\left(d_j - {\varepsilon \over 8} \right)
\right] F_0(k)
\nonumber
\\
&\le&\left(d_i + {\varepsilon \over 2} \right) F_0(k)
\le \left(1 + {\varepsilon \over 2} \right) CF(k),
\label{add-eqn-31}
\end{eqnarray}
where the last inequality is due to $d_j \le 1$ and $F_0(k) = C F(k)$
for all $k \ge k'$. Applying (\ref{add-eqn-31}) to (\ref{add-eqn-54}),
we have for $k \ge k'' > 2 k'$,
\begin{eqnarray}
F_i \ast F_j (k) 
&\le& \left(1 + {\varepsilon \over 2} \right) CF(k) 
+ \sum_{l=0}^{k'-1}  F_i(l)F_j(k - l) 
\nonumber
\\
&\le& \left(1 + {\varepsilon \over 2} \right) CF(k) 
+ \sup_{k - k'+1 \le l \le k}F_j (l).
\label{add-eqn-29}
\end{eqnarray}
Further for $k \ge k'' > 2k'$, $k - k'+1 > k'+1$ and thus
(\ref{eq-add-asum0}) and (\ref{eq-add-asum1}) yield
\begin{eqnarray}
\sup_{k - k'+1 \le l \le k}F_j (l)
&\le&
\left(1 + {\varepsilon \over 2} \right)\sup_{k - k'+1 \le l \le k}F_0 (l)
\nonumber
\\
&=&
\left(1 + {\varepsilon \over 2} \right)
\sup_{k - k'+1 \le l \le k}{F_0(l) \over F_0(k)} \cdot CF(k)
\nonumber
\\
&\le&
\left(1 + {\varepsilon \over 2} \right)
{1 \over (1 - \varepsilon)^{k'-1}} \cdot CF(k)
\nonumber
\\
&=& \left(1 + {\varepsilon \over 2} \right)C_{\varepsilon}' \cdot CF(k),
\quad k \ge k'' > 2k',
\label{add-eqn-28}
\end{eqnarray}
where $C_{\varepsilon}' = 1/(1-\varepsilon)^{k'-1}$.
Substituting (\ref{add-eqn-28}) into (\ref{add-eqn-29}), we obtain
\begin{eqnarray}
F_i \ast F_j (k)  
&\le& \left(1 + {\varepsilon \over 2} \right) 
\left( 1 + C_{\varepsilon}' \right) CF(k)
\nonumber
\\
&\le& ( 1 + \varepsilon ) 
\cdot 2C_{\varepsilon}' CF(k)
\nonumber
\\
&\le& 2 C_{\varepsilon}' \cdot (1 + \varepsilon)^2 CF(k),
\quad k \ge k'',
\label{add-eqn-25}
\end{eqnarray}
where we use $C_{\varepsilon}' \ge 1$. Note here that $F_i \ast F_j
(k) \le 1$ for all $k \in \bbZ_+$ and
\[
\sup_{k_0 \le k \le k''-1}F(k) / F(k'') \in (0,\infty).
\]
Therefore there exists some $C_{\varepsilon}'' > 0$ such that
\begin{eqnarray}
F_i \ast F_j (k) 
&\le& C_{\varepsilon}'' {F(k) \over F(k'')}
\nonumber
\\
&\le& {C_{\varepsilon}'' \over CF(k'')} \cdot (1 + \varepsilon)^2 CF(k),
~~~k_0 \le k \le k''-1. 
\label{add-eqn-26}
\end{eqnarray}
We now define $K_{\varepsilon}$ as
\[
K_{\varepsilon}
= \max\left(2 C_{\varepsilon}', 
{C_{\varepsilon}'' \over CF(k'')},
{2+\varepsilon \over \varepsilon(1+\varepsilon)^2}C_{\varepsilon}' 
\right).
\]
We then have the following inequality (which is used later).
\begin{equation}
\left( 1 + {\varepsilon \over 2} \right)  C_{\varepsilon}'
\le  K_{\varepsilon} (1 + \varepsilon)^2
{\varepsilon \over 2}.
\label{add-eqn-55}
\end{equation}
Further combining (\ref{add-eqn-25}) and (\ref{add-eqn-26}) leads to
\begin{equation}
F_i \ast F_j (k) 
\le K_{\varepsilon} (1 + \varepsilon)^2 CF(k),
\qquad k \ge k_0.
\label{add-eqn-27}
\end{equation}

Next we show (\ref{eq-prop-F(k)-bound}) for the convolution of three
mass functions $F_i, F_j$ and $F_{\nu}$ ($i,j,\nu=0,1,\dots,m$).  It
follows from (\ref{add-eqn-27}) and $F_0(k) = C F(k)$ for all $k \ge
k'$ that
\[
F_i \ast F_j (k) 
\le K_{\varepsilon} (1 + \varepsilon)^2 F_0(k),
\qquad k \ge k'.
\]
From this and (\ref{add-eqn-28}), we have for $k \ge k'' > 2k'$,
\begin{eqnarray}
\lefteqn{
F_i \ast F_j \ast F_{\nu} (k) 
}
\quad &&
\nonumber 
\\
&=&
\sum_{l=0}^{k - k'} F_i \ast F_j (k-l) F_{\nu}(l) 
+ \sum_{l=0}^{k'-1} F_i \ast F_j(l) F_{\nu}(k-l)
\nonumber 
\\
&\le&
\sum_{l=0}^{k - k'} F_i \ast F_j (k-l) F_{\nu}(l) 
+ \sup_{k - k'+1 \le l \le k} F_{\nu}(l)
\nonumber 
\\
&\le&
K_{\varepsilon} (1 + \varepsilon)^2 \sum_{l=0}^{k - k'} F_0(k-l) F_{\nu}(l) 
+ \left( 1 + {\varepsilon \over 2} \right)  C_{\varepsilon}' CF(k).
\label{add-eqn-30}
\end{eqnarray}
Applying (\ref{add-eqn-31}) and (\ref{add-eqn-55}) to
(\ref{add-eqn-30}) yields for $k \ge k'' > 2k'$,
\begin{eqnarray}
\lefteqn{
F_i \ast F_j \ast F_{\nu}(k) 
}
\quad &&
\nonumber 
\\
&\le&
K_{\varepsilon} (1 + \varepsilon)^2  
\left( 1 + {\varepsilon \over 2} \right) C F(k) 
+ K_{\varepsilon} ( 1 + \varepsilon )^2
{\varepsilon \over 2} CF(k)
\nonumber
\\
&=& 
K_{\varepsilon} (1 + \varepsilon)^2  
\left( 1 + {\varepsilon \over 2}  + {\varepsilon \over 2}\right) CF(k)
\nonumber
\\
&=& K_{\varepsilon} (1 + \varepsilon)^3 CF(k).
\nonumber
\end{eqnarray}
Further using $C_{\varepsilon}'' > 0$ such that (\ref{add-eqn-26}) holds, 
we obtain
\begin{eqnarray*}
F_i \ast F_j \ast F_{\nu}(k)
&\le& C_{\varepsilon}'' {F(k) \over F(k'')}
\nonumber
\\
&\le& {C_{\varepsilon}'' \over CF(k'')} \cdot (1 + \varepsilon)^3 CF(k),
~~~k_0 \le k \le k''-1. 
\end{eqnarray*}
Therefore $F_i \ast F_j \ast F_{\nu}(k) \le K_{\varepsilon} (1 +
\varepsilon)^3 CF(k)$ for $k \ge k_0$.

By repeating the above argument, we can prove that
(\ref{eq-prop-F(k)-bound}) holds under the additional condition that
$c_j > 0$ for all $j=1,2,\dots,m$.  In what follows, we remove this
condition.

Without loss of generality, we assume that $c_j = 0$ for
$j=1,2,\dots,m'$ ($1\le m' \le m$) and $c_j > 0$ for
$j=m'+1,m'+2,\dots,m$. Then for any $\delta > 0$, there exists some
positive integer $k_{\ast}:=k_{\ast}(\delta) \ge k_0$ such that for
all $k \ge k_{\ast}$,
\[
F_j(k) \le \delta F(k),\qquad j = 1,2,\dots,m'.
\]
Let $\{\tilde{F}_j(k); k\in \bbZ_+\}$ ($j = 1,2,\dots,m'$) denote a
probability mass function such that
\[
\tilde{F}_j(k)
=\left\{
\begin{array}{ll}
F_j(k)/\Theta_j, & k < k_{\ast},
\\
\delta F(k)/\Theta_j, & k \ge k_{\ast},
\end{array}
\right.
\]
where $\Theta_j:=\Theta_j(\delta) =
\sum_{k=0}^{k_{\ast}-1}F_j(k)+\sum_{k=k_{\ast}}^{\infty}\delta
F(k)$. It then follows that $F_j(k) \le \Theta_j\tilde{F}_j(k)$ for
all $k \in \bbZ_+$ and $j = 1,2,\dots,m'$. Thus we have
\begin{eqnarray}
\lefteqn{
F_1^{\ast n_1} \ast F_2^{\ast n_2} 
\ast \cdots \ast F_m^{\ast n_m}(k) 
}
\quad &&
\nonumber
\\
&\le& \prod_{j=1}^{m'}\Theta_j^{n_j} \cdot 
\tilde{F}_1^{\ast n_1} \ast \cdots \ast \tilde{F}_{m'}^{\ast n_{m'}} 
\ast F_{m'+1}^{\ast n_{m'+1}} \ast \cdots \ast F_m^{\ast n_m}(k).
\label{add-eqn-32}
\end{eqnarray}
By definition,
\[
\lim_{k\to\infty} {\tilde{F}_j(k) \over F(k)}
= {\delta \over \Theta_j} > 0,
\qquad j = 1,2,\dots,m'.
\]
Therefore for any $\varepsilon > 0$, there exists some
$C_{\varepsilon} > 0$ such that
\begin{eqnarray}
\lefteqn{
\tilde{F}_1^{\ast n_1} \ast \cdots \ast \tilde{F}_{m'}^{\ast n_{m'}} 
\ast F_{m'+1}^{\ast n_{m'+1}} \ast \cdots \ast F_m^{\ast n_m}(k)
}
\qquad \qquad\qquad\qquad &&
\nonumber
\\
&\le& 
C_{\varepsilon} (1+\varepsilon)^{n_1+n_2+\cdots+n_m}F(k).
\label{add-eqn-56}
\end{eqnarray}
Note here that $\lim_{\delta \downarrow 0}\Theta_j(\delta) = 1$ for all
$j=1,2,\dots,m'$.  Substituting (\ref{add-eqn-56}) into
(\ref{add-eqn-32}) and letting $\delta \downarrow 0$ yields
(\ref{eq-prop-F(k)-bound}).

\section{Examples}\label{appendix-example}

\subsection{M/GI/1 queue with Pareto service-time distribution}
\label{appendix-MG1}

We consider a stable M/GI/1 queue with a Pareto service-time
distribution. Let $\lambda$ denote the arrival rate of customers. Let
$H$ denote the service time distribution, which is given by
\[
H(x) = 1 - (x+1)^{-\gamma},
\quad x \ge 0,
\]
with $\gamma > 1$ and $\gamma \not\in \bbN$. Note here that the mean
service time is equal to $1/(\gamma - 1)$ and thus the traffic
intensity, denoted by $\rho$, is equal to $\lambda/(\gamma-1) < 1$.
Let $\widetilde{H}(s)$ denote the Laplace-Stieltjes transform (LST) of
the service time distribution $H$. It then follows from Theorem~8.1.6
in \cite{Bing89} that
\begin{equation}
\widetilde{H}(s) 
= \sum_{j=0}^{\lfloor \gamma \rfloor} h_j {(-s)^j \over j!} 
- \Gamma(1-\gamma)s^{\gamma} + o(s^{\gamma}),
\label{eqn-hat{H}(s)}
\end{equation}
where $h_j = \int_0^{\infty}x^j \rd H(x)$ $(j=1,2,\dots)$,
$f(x)=o(g(x))$ represents $\lim_{x\to0}f(x)/g(x)=0$ and $\Gamma$
denotes the Gamma function. Equation (\ref{eqn-hat{H}(s)}) yields
\begin{equation}
\widetilde{H}(\lambda - \lambda z) 
= \sum_{j=0}^{\lfloor \gamma \rfloor} h_j {(-\lambda)^j(1-z)^j \over j!} 
- \Gamma(1-\gamma)\lambda^{\gamma}(1-z)^{\gamma} + o((1-z)^{\gamma}).
\label{eqn-tilde{H}(lambda-lambdaz)}
\end{equation}
It is well-known that the stationary queue length distribution of the M/GI/1 queue,
denoted by $\{x(k);k\in \bbZ_+\}$, is identical with the stationary
distribution of the following stochastic matrix:
\[
\left(
\begin{array}{ccccc}
\alpha(0)    & \alpha(1)    & \alpha(2)    & \alpha(3)    & \cdots
\\
\alpha(0)    & \alpha(1)    & \alpha(2)    & \alpha(3)    & \cdots
\\
0      & \alpha(0)    & \alpha(1)    & \alpha(2)    & \cdots
\\
0      &   0    & \alpha(0)    & \alpha(1)    & \cdots
\\
\vdots & \vdots & \vdots & \vdots & \ddots
\end{array}
\right),
\]
where $\{\alpha(k);k\in\bbZ_+\}$ satisfies $\sum_{k=0}^{\infty}z^k
\alpha(k) = \widetilde{H}(\lambda - \lambda z)$ and thus
$\sum_{k=1}^{\infty}k\alpha(k)=~\rho$.

Let $\overline{\alpha}(k) = \sum_{l=k+1}^{\infty}\alpha_l$ for
$k\in\bbZ_+$. From (\ref{eqn-tilde{H}(lambda-lambdaz)}), we then have
\begin{eqnarray}
\sum_{k=0}^{\infty}z^k \overline{\alpha}(k) 
&=& {1 - \widetilde{H}(\lambda - \lambda z) \over 1 - z}
\nonumber
\\
&=& -\sum_{j=1}^{\lfloor \gamma \rfloor} h_j {(-\lambda)^j(1-z)^{j-1} \over j!} 
\nonumber
\\
&&{}
+ \Gamma(1-\gamma)\lambda^{\gamma}(1-z)^{\gamma-1} + o((1-z)^{\gamma-1}).
\label{add-eqn-22}
\end{eqnarray}
Applying Lemma~5.3.2 in \cite{Wilf94} to
(\ref{eqn-tilde{H}(lambda-lambdaz)}) and (\ref{add-eqn-22}) yields
\begin{eqnarray}
\alpha(k) 
&\simhm{k}& \gamma \lambda^{\gamma}  k^{-\gamma-1},
\label{asymp-alpha(k)}
\\
\overline{\alpha}(k) 
&\simhm{k}& \lambda^{\gamma}k^{-\gamma},
\label{asymp-bar-alpha(k)}
\end{eqnarray}
where $f(x) \simhm{x} g(x)$ represents $\lim_{x\to\infty}f(x)/g(x) =
1$.  Note that (\ref{asymp-alpha(k)}) shows that the discrete
distribution $\{\alpha(k);k\in\bbZ_+\}$ is in the class $\calL_{\rm
  loc}$. In fact, as shown later, $\{\alpha(k)\} \in \calS^{\ast}$,
i.e., $\{\alpha_{\rm e}(k)\} \in \calS_{\rm loc}(1)$, where
$\alpha_{\rm e}(k) = \overline{\alpha}(k) / \rho$ for
$k=0,1,\dots$. Therefore it follows from Theorem~\ref{sec4-thm-sub}
that
\[
x(k) 
\simhm{k} {\rho \over 1 - \rho} \cdot \alpha_{\rm e}(k)
= {\rho \over 1 - \rho} \cdot {\overline{\alpha}(k) \over \rho}
\simhm{k} {\lambda^{\gamma} \over 1-\rho} k^{-\gamma}.
\]

In what follows, we prove that $\{\alpha(k)\} \in \calS^{\ast}$, i.e., 
\[
\sum_{l=0}^k \overline{\alpha}(l) \overline{\alpha}(k-l)
 \simhm{k} 2\rho\cdot \overline{\alpha}(k).
\]
Let $\nu:=\nu(k)$ denote an integer such that $k/3 \le \nu(k) < k/2$.
For $k \in \bbZ_+$, we have
\begin{eqnarray}
\sum_{l=0}^k {\overline{\alpha}(l) \overline{\alpha}(k-l) \over \overline{\alpha}(k)}
= 2\sum_{l=0}^{\nu-1} 
\overline{\alpha}(l){\overline{\alpha}(k-l) \over \overline{\alpha}(k)} 
+ \sum_{l=\nu}^{k-\nu} 
\overline{\alpha}(l){\overline{\alpha}(k-l) \over \overline{\alpha}(k)}. 
\label{eqn-00}
\end{eqnarray}
From (\ref{asymp-bar-alpha(k)}), we obtain 
\begin{eqnarray}
\lim_{k\to\infty}\sum_{l=0}^{\nu-1} 
\overline{\alpha}(l){\overline{\alpha}(k-l) \over \overline{\alpha}(k)}
&=& \sum_{l=0}^{\nu-1} \overline{\alpha}(l) 
\lim_{k\to\infty}{\overline{\alpha}(k-l) \over \overline{\alpha}(k)}
=\sum_{l=0}^{\nu-1} \overline{\alpha}(l).
\label{eqn-01}
\end{eqnarray}
Further it follows from (\ref{asymp-bar-alpha(k)}) that for any
$\varepsilon > 0$ there exists some $k_{\ast} \in \bbZ_+$ such that
for all $k \ge k_{\ast}/3$,
\[
1 - \varepsilon 
< {\overline{\alpha}(k) \over \lambda^{\gamma}k^{-\gamma}} < 1 + \varepsilon,
\]
which implies that for $k \ge k_{\ast}$ and $k/3 \le \nu < k/2$,
\begin{eqnarray}
\sum_{l=\nu}^{k-\nu} \overline{\alpha}(l)
{\overline{\alpha}(k-l) \over \overline{\alpha}(k)}
&\le& {(1+\varepsilon)^2 \over 1-\varepsilon}
\sum_{l=\nu}^{k-\nu}\lambda^{\gamma}l^{-\gamma}\left({k-l \over k}\right)^{-\gamma}
\nonumber
\\
&\le& {(1+\varepsilon)^2 \over 1-\varepsilon}
\lambda^{\gamma}(k-2\nu+1)\nu^{-\gamma}
\left({\nu \over k}\right)^{-\gamma}
\nonumber
\\
&\le& {(1+\varepsilon)^2 \over 1-\varepsilon}
\lambda^{\gamma} k \left({k \over 3}\right)^{-\gamma} 3^{\gamma}
\nonumber
\\
&\le& {(1+\varepsilon)^2 \over 1-\varepsilon}
(9\lambda)^{\gamma} k^{-\gamma+1} \to 0,\quad \mbox{as}~k\to\infty.
\label{eqn-02}
\end{eqnarray}
Finally, applying (\ref{eqn-01}) and (\ref{eqn-02}) to (\ref{eqn-00})
and letting $\nu \to \infty$ yield
\[
\lim_{k\to\infty}
\sum_{l=0}^k 
{\overline{\alpha}(l) \overline{\alpha}(k-l) \over \overline{\alpha}(k)}
= 2\sum_{l=0}^{\infty} \overline{\alpha}(l) = 2\rho.
\]

\subsection{Discrete-time queue with disasters and Pareto-distributed batch arrivals}

This subsection considers a discrete-time single-server queue with
disasters and Pareto-distributed batch arrivals. The time interval
$[n,n+1)$ ($n\in\bbZ_+$) is called slot $n$. Customers and disasters
  can arrive at the beginnings of respective slots, whereas departures
  of served customers can occur at the ends of respective slots.

We assume that the numbers of customer arrivals in respective slots
are independent and identically distributed (i.i.d.) with a discrete
Pareto distribution, $\beta(k) = 1/(k+1)^{\gamma} - 1/(k+2)^{\gamma}$
($k\in\bbZ_+$), where $\gamma > 1$. Service times are i.i.d. with a
geometric distribution with mean $1/(1-q)$ ($0 < q < 1$).  We also
assume that at most one disaster occurs at one slot with probability
$\phi$ ($0 < \phi < 1$), which are independent of the arrival process
of customers. If a disaster occurs in a slot, then both customers
arriving in the slot and all the ones in the system are removed.

Let $L_n$ ($n \in \bbZ_+$) denote the number of customers at the
middle of slot $n$. It then follows from
Proposition~\ref{prop-stability-2} that $\{L_n;n\in\bbZ_+\}$ is an
ergodic Markov chain whose transition probability matrix is given by
\[
\left(
\begin{array}{cccccc}
b(0) 			& b(1) & b(2) & b(3) & b(4) & \cdots
\\
\phi + a(0)  	& a(1) & a(2) & a(3) & a(4) & \cdots
\\
\phi    		& a(0) & a(1) & a(2) & a(3) & \cdots
\\
\phi    		& 0    & a(0) & a(1) & a(2) & \cdots
\\
\phi    		& 0    & 0	  & a(0) & a(1) & \ddots
\\
\vdots   		& \vdots   & \vdots  & \vdots  & \ddots & \ddots 
\end{array}
\right),
\]
where
\begin{align*}
b(0) &= \phi + (1 - \phi)\beta(0), 
\\
b(k) &= (1 - \phi)\beta(k), & k &=1,2,\dots, 
\\
a(0) &= (1 - \phi)\beta(0)(1-q), 
\\
a(k) &= (1 - \phi)[\beta(k-1)q + \beta(k)(1-q)], & k &=1,2,\dots.
\end{align*}
It is easy to see that $\sum_{k=0}^{\infty}a(k)= 1 - \phi
< 1$ and
\[
\lim_{k\to\infty}{a(k) \over \beta(k)} 
= 1 - \phi,
\quad
\lim_{k\to\infty}{b(k) \over \beta(k)} 
=  1-\phi.
\]
Note here that $\{\beta(k);k\in\bbZ_+\}$ is decreasing and
\[
\beta(k) \simhm{k} \gamma k^{-\gamma -1}.
\]
Thus as in subsection~\ref{appendix-MG1}, we can show that
$\{\beta(k);k\in\bbZ_+\} \in \calS_{\rm loc}(1)$. As a result, 
Theorem~\ref{thm-loc-substoch}~yields
\[
\lim_{n\to\infty}\PP(L_n=k)
\simhm{k} {1 - \phi \over \phi}\beta(k)
\simhm{k} {1 - \phi \over \phi}\gamma k^{-\gamma -1}.
\]

\section*{Acknowledgment}
The authors thank Professors Bara Kim and Jeongsim Kim for providing a
copy of their manuscript prior to publication.  Research of the second
author was supported in part by Grant-in-Aid for Young Scientists (B)
of Japan Society for the Promotion of Science under Grant
No.~24710165.

%
%
%

%
%

%

\end{document}